\documentclass[aop,preprint]{imsart}

\usepackage{amsfonts,amsmath,latexsym,amssymb}    
\usepackage{graphicx}
\RequirePackage[OT1]{fontenc}
\RequirePackage{amsthm,amsmath}
\RequirePackage[numbers]{natbib}
\RequirePackage[colorlinks,citecolor=blue,urlcolor=blue]{hyperref}

\arxiv{arXiv:1404.5705v7}

\startlocaldefs
\numberwithin{equation}{section}
\theoremstyle{plain}
\newtheorem{theorem}{Theorem}[section]
\newtheorem{lemma}[theorem]{Lemma}
\newtheorem{proposition}[theorem]{Proposition}

\newtheorem{remark}[theorem]{Remark}

\def\definedas{\stackrel{\Delta}{=}}

\newcommand{\real}{\mathbb{R}}

\newcommand{\abbr}[1]{{\sc\lowercase{#1}}}

\def\complex{\mathop{\raise .45ex\hbox{${\bf\scriptstyle{|}}$}
     \kern -0.40em {\rm \textstyle{C}}}\nolimits}

\newcommand{\calC}{{\cal C}}

\newcommand{\calF}{{\cal F}}
\newcommand{\calG}{{\cal G}}
\newcommand{\calH}{{\cal H}}

\newcommand{\calL}{{\cal L}}

\newcommand{\calN}{{\cal N}}
\newcommand{\calP}{{\cal P}}

\newcommand{\calU}{{\cal U}}

\newcommand{\calZ}{{\cal Z}}

\def\be{\beta}
\def\la{\lambda}

\numberwithin{equation}{section}

\def\definedas{\stackrel{\Delta}{=}}

\def\P{{\mathbb{P}}}
\def\E{{\mathbb{E}}}

\def\c{\theta}
\def\cn{\c_N}
\def\S{{\mathbb{S}}}

\endlocaldefs
\begin{document}

\begin{frontmatter}
\title{Component sizes for large quantum Erd\H{o}s-R\'{e}nyi graph near criticality}
\runtitle{Critical quantum random graphs}

\begin{aug}
\author{\fnms{Amir} \snm{Dembo}\thanksref{t1, m1}\ead[label=e1]{amir@math.stanford.edu}},
\author{\fnms{Anna} \snm{Levit}\thanksref{m2}\ead[label=e2]{anna.levit@math.ubc.ca}}
\and
\author{\fnms{Sreekar} \snm{Vadlamani}\thanksref{m3}
\ead[label=e3]{sreekar@tifrbng.res.in}
}

\thankstext{t1}{Research funded by NSF grants DMS-11-06627 and DMS-1613091}
\runauthor{A. Dembo et al.}

\affiliation{Stanford University\thanksmark{m1}, University of British Columbia\thanksmark{m2} and TIFR-CAM\thanksmark{m3}}

\address{Department of Statistics \\
\ \ \ and Department of Mathematics,\\
Stanford University, Stanford, CA 94305\\
\printead{e1}}

\address{Department of Mathematics,\\
The University of British Columbia,\\
\#121-1984 Mathematics Road,\\
Vancouver, BC, V6T 1Z2, Canada\\
\printead{e2}}

\address{TIFR-Center for Applicable Mathematics, \\ Post Bag 6503, GKVK Post Office,\\
Bangalore 560065\\
\printead{e3}
}
\end{aug}

\begin{abstract}
The $N$ vertices of a quantum random graph are each a circle independently 
punctured at Poisson points of arrivals, with parallel connections derived through 
for each pair of these punctured circles by yet another 
independent Poisson process. Considering these 
graphs at their
critical parameters, we show that 
the joint law of the re-scaled by $N^{2/3}$ and ordered sizes of their 
connected components,  
converges to that of the ordered lengths of excursions above 
zero for a reflected Brownian motion with drift. 
Thereby, this work forms the first example of an inhomogeneous random graph, beyond the 
case of effectively rank-1 models, which is rigorously shown to be in the 
Erd\H{o}s-R\'{e}nyi graphs universality class in terms of Aldous's results.
\end{abstract}



\end{frontmatter}

\section{Introduction}
The Erd\H{o}s-R\'{e}nyi random graph \cite{erdos-renyi} is the simplest and
most studied example of a random graph ensemble. 
Such a graph, denoted by $G(N,p)$, has $N$ vertices, with each pair of  
vertices connected with probability $p$, independently of all other pairs. 
Its phase transition phenomena are well understood. In particular, for 
$p=\frac{c}{N}$ with $c>1$, 
the largest component in $G(N,p)$ has $\Theta(N)$ vertices and 
the second largest ${\rm O}(\ln N)$ vertices (as $N\rightarrow\infty$, with probability $1$), 
for $p=\frac{c}{N}$ with $c<1$ the largest component has ${\rm O}(\ln N)$ vertices and
when $p=\frac{1}{N}$, the largest component of $G(N,p)$ has 
$\Theta(N^{2/3})$ vertices (c.f. \cite{erdos-renyi, bollobas, luczak}).

Aldous \cite{aldous} considered the asymptotic behavior of $G(N,p)$ 
inside the ``scaling window" of this phase transition, namely for $N \to \infty$
and $|p-1/N|$ small enough, showing that the ordered set of 
component sizes rescaled by $N^{2/3}$ then converges to an ordered set of excursion
lengths of reflected inhomogeneous Brownian motion with a certain drift.
Various other random graph models exhibit a phase transition 
phenomenon similar to the Erd\H{o}s-R\'{e}nyi random graph. While some 
further follow the same behaviour as $G(N,p)$ in their near-critical regime,
the near-critical regime of others falls 
into different universality classes.     

For example, Nachmias and Peres \cite{nach-peres2} prove that the random graph ensemble obtained 
by performing percolation on a random $d$-regular ($d\geq 3$) graph on $N$ vertices with percolation probability $p=1/(d-1) + a N^{-1/3}$ 
for $a\in \mathbb{R}$ fixed, falls 
into the same universality class as the Erd\H{o}s-R\'{e}nyi random graph.
The random multi-graph whose $N$ vertices are constructed using the configuration model, with its vertex degrees being 
i.i.d. variables, each having the distribution $\nu$, has a richer behavior. 
Indeed, Joseph \cite{joseph} shows 
that when $\nu$ has a finite third moment,
the near-critical regime of this model falls into the Erd\H{o}s-R\'{e}nyi random graph's universality class, whereas  
if $\nu_k\sim c k^{-\tau}$ as $k\rightarrow \infty$, $c>0$, $\tau\in(3,4)$ then the 
relevant scaling changes to $N^{-(\tau-2)/(\tau-1)}$ and the limit is 
an ordered set of the excursion lengths of some other drifted process with independent increments 
above past minima 
(near-critical regime of the Erd\H{o}s-R\'{e}nyi universality class, is
also obtained in \cite{riordan} for more general class of 
degree distributions of finite third moment). 
A similar behavior has been found in the near critical regime of 
the Rank-1 model (a special case of the general in-homogeneous random graph 
studied in \cite{bollobas-janson-riordan}, which has received much attention 
recently). Such graph has random i.i.d. weights $\{x_i\}$ associated to its vertices, 
and edges chosen independently, with the edge $(i,j)$ chosen with probability
$p_{i,j}=\min\{c\frac{x_i x_j}{N}, 1\}$, for some positive constant $c=c(N)$ (c.f. 
\cite{turova}). For $x_i$ having finite third moment, the near-critical regime 
corresponds to $c(N)=1+a N^{-1/3}$, in which case \cite{aldous} shows 
that this model (formulated slightly differently), falls into the 
Erd\H{o}s-R\'{e}nyi graph's universality class (similar results have been 
later proved in \cite{bhamidi-hofstad-leeuwaarden1,turova2}).
In contrast, for Rank-1 model with power-law degrees of exponent $\tau \in (3,4)$,
\cite{bhamidi-hofstad-leeuwaarden2} show that
the sizes of the components, re-scaled by $N^{-(\tau-2)/(\tau-1)}$, converge to 
hitting times of certain thinned L\'{e}vy process. 

Our aim here is to study the near-critical behavior of the 
so-called quantum version of Erd\H{o}s-R\'enyi random graph (\abbr{qrg}). We note in passing, that 
both the motivation and terminology come from the stochastic 
geometric (Fortuin--Kasteleyn type) representation of the 
quantum Curie-Weiss model at inverse temperature 
$\beta>0$ (we exclude here the ground state case of $\beta=\infty$), 
in transverse magnetic field of strength $\lambda>0$ 
(at $\lambda=0$ it reduces to the Erd\H{o}s-R\'{e}nyi ensemble, see
Remark \ref{rmk:er-la-zero}). 
We refer the reader to \cite{ioffe} for more information on 
such stochastic geometric representations (that were 
originally developed in \cite{aizenman-klein-newman,campanino-klein-perez} 
for the general ferromagnetic context), moving on instead, 
to the precise description of the \abbr{qrg} (as in \cite{ioffe-levit}).

\begin{figure}
\begin{center}
\includegraphics[width=0.7\textwidth]{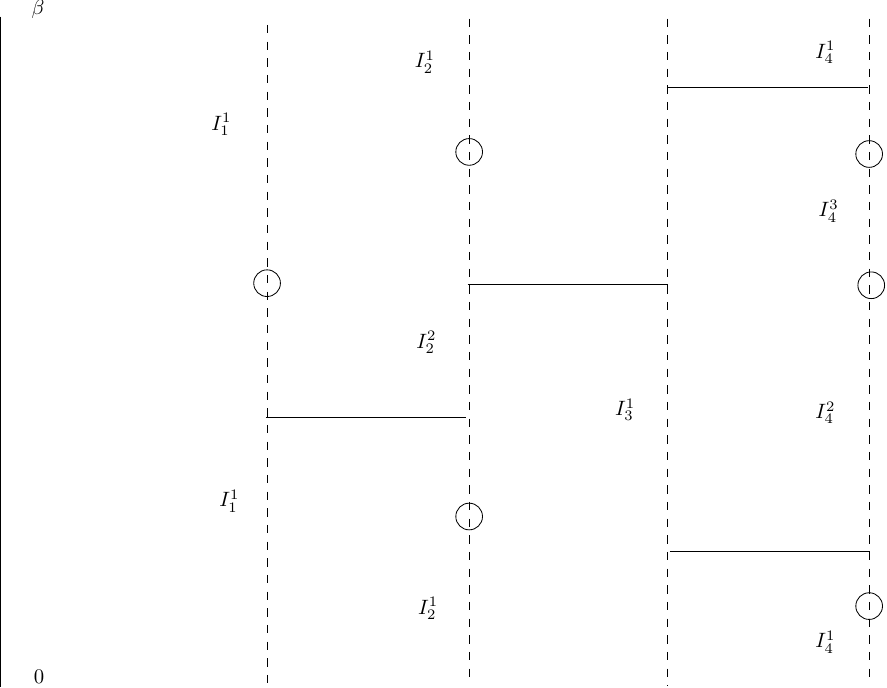}
\end{center}
\caption{ An example of the decomposition of $\calG_N^\beta$ after all the holes are punched and the links are drawn: 
$\calG_N^\beta\setminus\calH=\calC_1\vee\calC_2\vee\calC_3$, where 
$\calC_1=I_1^1\cup I_2^2\cup I_3^1\cup I_4^1\cup I_4^2$, $\calC_2=I_2^1$ 
and $\calC_3=I_4^3$.}  
\label{fig:1}
\end{figure}

\paragraph{The model:}
With $G_N = \{1,\ldots,N\}$ and $\S_{\be}$ denoting the circle of length $\beta$, 
let $\calG_N^\beta = G_N \times\S_{\be}$, 
associating to each site $i \in G_N$ the copy $\S^i_{\be}= i \times \S_{\be}$ 
of $\S_{\be}$, so a point in $\calG_N^\beta$ has two coordinates, its site 
(in $G_N$) and time (in $\S_\beta$) coordinates. 
The \abbr{qrg} is then the following random subset $\calG_N^\beta \setminus \calH$ 
of $\calG_N^\beta$, equipped with random links $\bigcup_{i,j} \calL_{i,j}$ 
between pairs of points of the type $\{(i,t)$ and $(j,t)$, for $i \ne j \}$.  
To construct the \abbr{qrg} we first punch within each $\S^i_{\be}$ 
finitely many holes, according to independent Poisson point processes $\calH_i$,
$i \in G_N$, of intensity $\lambda>0$, so
each resulting punctured circle $\S^i_{\be}\setminus\calH_i$ consists of $m_i$ 
disjoint connected intervals
\begin{equation}
\label{holedecomp}
\S^i_{\be}\setminus\calH_i=\bigcup_{l=1}^{m_i} I^l_i
\end{equation}
(the number of holes $\#\calH_i=m_i$, except when $\#\calH_i=0$, in which case
$m_i=1$). 
We next add links between pairs of points in $\calG_N^\beta$ 
of the same time coordinates (i.e. between points $(i,t)$ and $(j,t)$ where 
$i\neq j$ and $t\in \S_{\be}$), as follows. With each (unordered) pair of sites 
$i,j \in G_N$ we associate a copy $\S^{i,j}_{\be}$ of $\S_{\be}$ and a 
Poisson point process of links $\calL_{i,j}$ on $\S^{i,j}_{\be}$ with intensity $\frac1N$. 
The processes $\calL_{i,j}=\calL_{j,i}$ are assumed to be independent 
for different $(i, j)$ and also independent of the processes of holes $\calH_i$. 
Two intervals $I^l_i$ and $I^k_j$ of the decomposition \eqref{holedecomp} are then 
considered to be directly connected if there exists some 
$t \in \calL_{i,j}$ such that both $(i, t)\in I^l_i$ and $(j, t)\in I^k_j$. 
Setting $\calH := \cup_i \calH_i$ (a \emph{finite} collection of points),
the decomposition 
\begin{equation}
\label{decomp}
\calG_N^\beta\setminus\calH=\calC_1\vee \cdots \vee\calC_\ell
\end{equation}
of $\calG_N^\beta\setminus\calH$ into maximal connected components is, thereby, well 
defined
(see Figure~\ref{fig:1} for an example with $N = 4$). Further, each fixed 
$x \in\calG_N^\beta$ is a.s. not in $\calH$, 
hence the notion of the connected component $\calC(x)$ containing $x$ in 
the decomposition \eqref{decomp}, is also well defined, and hereafter the
\emph{size} of a connected component $\calC_j$ (or $\calC(x)$), means
the number of intervals it contains, and $\calP(\calC(x)) = \sum_I |I| \, 1_{I \in \calC(x)}$ denotes the
cumulative length of intervals constituting the component $\calC(x)$.

\begin{remark}\label{rmk:er-la-zero}
For $\lambda = 0$ there are no holes, so each $\S^i_{\be}\setminus\calH_i$ consists
of one connected component, which equals to $\S^i_{\be}$ itself. 
We are then back to the Erd\H{o}s-R\'{e}nyi random graph $G(N,p)$ with 
$p=1- e^{-\frac{\beta}{N}}$  (the probability that 
$\S^i_{\be}$ and $\S^j_{\be}$ are directly connected). 
\end{remark}

Treating each interval $I_i^k$ as a vertex, Janson in \cite{janson} notices that the 
\abbr{qrg} is an instance of the general in-homogeneous model 
of \cite{bollobas-janson-riordan}. However, the probability of direct 
connection between two intervals depends on the size of their overlap
and not only on the individual lengths of these intervals. Beyond
separating our model from the class of rank-1 models, this property 
makes it inherently different from the other models we have mentioned 
thus far (all of whom mimic
the idea of rank-1 random graphs, in the sense that certain 
\emph{vertex related weights} determine the probabilities in which edges are 
present in the graph).

An equivalent description of the \abbr{qrg} in case $\lambda>0$, 
which we adopt hereafter, has $N$ circles of length 
$\c\definedas\lambda\beta$ with a unit intensity Poisson 
process of holes on them, using now i.i.d. Poisson processes of intensity
$1/(\lambda N)$ for creating links between each pair of (punched) circles.
The critical curve for the \abbr{qrg} model
in the $(\beta,\lambda)$-parameter space, is obtained in \cite{ioffe-levit}
by comparisons with a critical branching process whose offspring distribution
is the \emph{cut-gamma} distribution $\Gamma_{\c}(2,1)$
(namely, the law of $J:=(J_- + J_+) \wedge \c$ for $J_-,J_+$ i.i.d. Exp($1$) variables). 
Using the preceding parametrization, the resulting curve 
$\beta=\beta_c(\lambda)$ corresponds to
\begin{equation}\label{eq:critical}
\beta_c =\frac{\c}{F(\c)}, \quad \lambda = F(\c), \quad \text{for} \quad  
F(\c) = 2 (1 - e^{-\c}) - \c e^{-\c} 
\end{equation}
(where $F(\c)$ is precisely the expected length $J$ 
of the interval $I$ in the
\abbr{QRG} upon our rescaling by $\lambda$).
It is easy to check that $\lambda (\c): [0,\infty) \mapsto [0,2)$ is concave, 
increasing and $\beta_c (\c) : [0,\infty) \mapsto [1,\infty)$ 
is strictly increasing, such that the curve 
$\beta_c (\lambda) : [0,2) \mapsto [1,\infty)$ is strictly increasing.
The critical curve is alternatively given by
\begin{equation}\label{eqn:crit:curve}
F(\beta,\lambda):= \lambda^{-1} F(\lambda \beta) = 1 \,,
\end{equation}
and it is further shown in \cite{ioffe-levit} that taking $F(\beta,\lambda)>1$ 
(equivalently, $\beta>\beta_c(\lambda)$), yields 
the emergence of an $\Theta(N)$-giant connected component in the 
disjoint decomposition \eqref{decomp}, whereas when
$F(\beta,\lambda)<1$ (equivalently, $\beta<\beta_c(\lambda)$),
all connected components are typically
of order ${\rm O}(\ln N)$. Our first result complements \cite{ioffe-levit} by 
proving that at criticality the largest component is of size $\Theta(N^{2/3})$
(so the \abbr{qrg} admits a version of the Erd\H{o}s-R\'{e}nyi phase transition).
\begin{theorem} \label{thm:crit:size}
Suppose $(\beta,\lambda)$ is a critical point, namely $F(\beta,\lambda)=1$. 
Then, for the largest component 
$\calC_{\rm max}$ of the \abbr{qrg}, we have that: 
\begin{itemize}
\item[(a)] There exist $c_*$, $N_0$ and $A_0$ finite, such that for all 
$N>N_0$ and $A>A_0$, 
\begin{equation}\label{eq:pmax-ubd}
\P\left( \calP(\calC_{\rm max}) > A N^{2/3}\right) \le c_* A^{-3/2}.
\end{equation}
\item[(b)] There exists $N_1$ finite such that for all 
$N>N_1$ and $\delta>0$,
\begin{equation}\label{eq:pmax-lbd}
\P\left( \calP(\calC_{\rm max}) <  \lfloor\delta N^{2/3}\rfloor \right) \le 
(6+4\beta^2) \delta^{3/5}.
\end{equation}
\end{itemize}
\end{theorem}

Our primary objective is to further 
analyze the \abbr{qrg} model, 
and in particular its component sizes near criticality, 
thereby confirming that the \abbr{qrg} is in the same universality class as the 
Erd\H{o}s-R\'{e}nyi random graph. Whereas our proofs also rely on
an exploration process for estimating the connected components sizes,
in contrast to all cases dealt with before (i.e. 
\cite{aldous, bhamidi-hofstad-leeuwaarden1, joseph, nach-peres1, nach-peres2}), 
here we may have many intervals sharing the same vertex 
(that is, a site $i \in G_N$, or alternatively, 
the corresponding circle $\S_\c^i$). Thus, our 
exploration process (or breadth first walk), may
re-visit an already visited vertex (circle), as many 
times as the number of intervals sharing such vertex. The 
latter is an \emph{unbounded} random variable, thereby posing a 
serious challenge to our analysis. While Theorem \ref{thm:crit:size} 
is rough enough that we can surpass this difficulty by showing that 
multiple returns to same vertex are rare enough to not matter, 
this is no longer true for our main result, Theorem \ref{thm:joint:cgs}, 
about the scaling limits of ordered component sizes. Indeed, our limiting process drift differs from that
of \cite{aldous} by additional quadratic factor 
representing the already explored portion 
of the relevant circle.
Indeed, the question of convergence of such quantum random graphs, as metric spaces, 
is completely open due to this precise problem of multiple visits to the same vertex.

The following quantities are required for our main result. First, let 
\begin{align}\label{eqn:sig-dfn}
\sigma^2(\c) & =  \frac{\E [ J^2 ]}{F(\c)^2} \,, \\
\gamma(\c) &= \frac{\E[\widehat{F}(\c-J)]}{\c F(\c)} \,,
\quad \text{ for } \quad 
\widehat{F}(x)=2(x-1+e^{-x}) + \frac{x^3}{3} e^{-\c}\,.
\label{eqn:gamma-dfn}
\end{align}
Then, for standard Brownian motion $\{W(s), s \ge 0 \}$ and any $a \in \real$, 
consider the processes
\begin{align}\label{eqn:def:W-a-beta}
W^{a,\c}(s) &:= \sigma(\c) W(s) + \rho^{a,\c}(s) \,, \\
\rho^{a,\c}(s) & := a s -\frac{s^2}{2} \Big( 1 - \gamma(\c)
 \Big) \,,
\label{eqn:def:rho-a-beta}
\end{align}
and the associated process of non-negative excursions
\begin{equation}\label{eqn:def:B-a-beta}
B^{a,\c}(s) = W^{a,\c}(s) - \min_{0\le u\le s}W^{a,\c}(u).
\end{equation}

\begin{theorem}\label{thm:joint:cgs}
Fix $a\in \mathbb{R}$ and $(\beta,\lambda)$ a point on the critical curve
of \eqref{eqn:crit:curve}. Consider parameters
$(\beta_N,\lambda_N) \to (\beta,\lambda)$ such that 
$F(\beta_N,\lambda_N) = 1 + a N^{-1/3}$. Then, denoting 
the ordered sizes of components 
of the graph by $|\calC^{a,N}_1|,\,|\calC^{a,N}_2|,\ldots$, 
we have when $N\rightarrow \infty$ that
\[
\left(N^{-2/3}|\calC^{a,N}_1|,\, N^{-2/3}|\calC^{a,N}_2|,\ldots\right)\,\,\, \stackrel{d}{\Rightarrow} \,\,\, (\gamma_1,\gamma_2,\ldots)\,,
\]
where $\{\gamma_j\}$ denote the ordered lengths of the excursions of the process
$B^{a,\c}$ above zero, and the convergence of component sizes holds 
with respect to the $l^2_{\searrow}$ topology (as defined in \cite{aldous}).
\end{theorem}

In Section \ref{sec:crit:scale} we prove Theorem \ref{thm:crit:size} by adapting 
to our context the ideas set forth in \cite{nach-peres1}. Specifically, the main 
task here is to construct a pair of manageable auxiliary counting processes,
which are not too far apart, while stochastically dominating (from above and below, 
respectively), the counting process that determines the size of our components (thereby 
circumventing much of the difficulty associated with the precise counting).

Section \ref{sec:joint:cgs} is devoted to the proof of Theorem \ref{thm:joint:cgs}
which requires finer estimates and thereby some new ideas. What sets  
our analysis apart of all those mentioned before, is Proposition
\ref{prop:nu-m} which provides rough estimates on the number of sites of \abbr{qrg}
visited twice, or
more,
during the first $k = {\rm O}(N^{2/3})$ steps of the exploration process. 
It shows in particular that only the \emph{first return to a site} 
plays a crucial role, with subsequent returns playing no role when 
the relevant limit is considered. Combined with further rough estimate
on the number of sites visited 
exactly once during the first $k$ steps, it
thus allows us to thereafter adapt the program of \cite{aldous}
to the \abbr{qrg} setting. Specifically, 
Subsection \ref{subsec:brownian:excursions} deals with
weak convergence of the law induced by the rescaled breadth first walk 
to the law of $W^{a,\c}$ defined on the space
of \abbr{RCLL} 
functions $D\left([0,\infty)\right)$, equipped with the topology of 
uniform convergence on finite intervals. Finally,
in Subsection \ref{subsec:joint:convergence} we collate all the above results into
a proof of Theorem \ref{thm:joint:cgs}.

A further insight gained from our proof is that the \abbr{qrg} model
is in the Erd\"os-R\'enyi universality class by the confluence of
two reasons: first the small probability of many returns to the 
same vertex (circle); second and more crucial is the relatively fast relaxation of 
its exploration process, which thereby behaves approximately as a Markov process. 
One may examine the latter feature in many other inhomogeneous 
random graph models, and where it is present, proceed to try proving that
they too belong to the Erd\"os-R\'enyi universality class.

\paragraph{Acknowledgment} 
This work benefited from helpful discussions of A.L. with Omer Angel 
and from the hosting of A. L. and A. D. by the
Mathematical Sciences Research Institute (as part of its program 
on random spatial processes). The authors would also like to thank
J\'ulia Komj\'athy for useful comments.


\section{Proof of Theorem \ref{thm:crit:size}}
\label{sec:crit:scale}

We shall first prove that $|\calC_{\max}| = \Theta_{\P}(N^{2/3})$, then to 
conclude the result of Theorem \ref{thm:crit:size}, we shall
use natural bounds arising from the arguments used to prove the former.

In particular, our first step towards proving Theorem \ref{thm:crit:size} will be
the following proposition.

\begin{proposition} \label{prop:crit:size}
Suppose $(\beta,\lambda)$ is a critical point, namely $F(\beta,\lambda)=1$. 
Then, for the largest component 
$\calC_{\rm max}$ of the \abbr{qrg}, we have that: 
\begin{itemize}
\item[(a)] There exist $c_*$, $N_0$ and $A_0$ finite, such that for all 
$N>N_0$ and $A>A_0$, 
\begin{equation}\label{eq:cmax-ubd}
\P\left( |\calC_{\rm max}| > A N^{2/3}\right) \le c_* A^{-3/2}.
\end{equation}
\item[(b)] There exists $N_1$ finite such that for all 
$N>N_1$ and $\delta>0$,
\begin{equation}\label{eq:cmax-lbd}
\P\left( |\calC_{\rm max}| <  \lfloor\delta N^{2/3}\rfloor \right) \le 
(6+4\beta^2) \delta^{3/5}.
\end{equation}
\end{itemize}
\end{proposition}

In proving $|\calC_{\max}| = \Theta_{\P}(N^{2/3})$, to bypass the problem of multiple visits 
of the same vertex by our exploration process
(as described in Subsection \ref{subsec:exploration}), we 
stochastically sandwich it between the over-counting process
of Subsection \ref{subsec:overcounting}, and the under-counting process
of Subsection \ref{subsec:free-exploration}. By combining the upper and 
lower bounds provided by these two auxiliary processes, we complete
the proof of Theorem \ref{thm:crit:size}.

\subsection{Exploration process}
\label{subsec:exploration}

Taking advantage of conditional independence properties of Poisson processes,
we start with an algorithmic definition of the exploration process for 
our model (following \cite{ioffe-levit}, who used it for examining 
a single component). This algorithm allows us to sequentially 
construct (or sample), the rescaled \abbr{qrg}, interval by interval.
In this description, our vertices (circles) have first been labeled 
$\{1,2\ldots,N\}$, and after $k$ steps of the algorithm, we fully 
explore $k$ intervals, having $A_k \ge 0$ active points
(unexplored ends of connections with the already explored intervals, or 
the point around which a new component starts),
while the rest of the space 
is declared to be neutral (note 
that we start exploring a new component of the graph
upon arriving at $A_{k-1}=0$, but not before).

\emph{\textbf{Initial stage:}} We fix the vertex $w_0 = 1$ and choose a point $t$ 
uniformly at random on this vertex. At end of step $k=0$ we have 
$A_0=1$, with $(w_0,t)$ as our sole active point and the whole
space considered neutral.

\emph{\textbf{At step $k \ge 1$:}} 
\begin{itemize}
\item[(a)] If $A_{k-1}>0$ we choose an active point $(w_k,t)$ whose vertex has the smallest index among all active points. 
In case of a tie, choose the active point which chronologically appeared earlier than the others 
on the same vertex. 
\item[(b)] If $A_{k-1}=0$ and there exists at least one neutral circle, we choose 
$w_k$ to be the neutral vertex with the smallest index and uniformly at random 
mark a new active point $(w_k,t)$ on this vertex.
\item[(c)] If $A_{k-1}=0$ and there is no neutral circle, we choose $w_k$ 
to be the vertex of smallest index among the vertices having some neutral part,
marking new active point $(w_k,t)$ uniformly at random on the 
neutral part of $\S_{\c}^{w_k}$.
\item[(d)] If $A_{k-1}=0$ and there is no neutral part available on any circle,
then this ends the exploration process.
\end{itemize}

Using i.i.d. Exp($1$) variables $J_-,J_+$, 
we carve out of the maximal neutral interval $\{w_k\} \times (t_1,t_2)$ 
around $(w_k,t)$, the sub-interval $I_k:=\{w_k\} \times \widetilde{I}$ 
for $\widetilde{I} = (t_1 \vee (t- J_-), t_2 \wedge (t + J_+))$. For 
$\S_{\c}^{w_k}$ completely neutral (apart from active points),
we take $t_2=-t_1=\infty$ and
$\widetilde{I} = \S_{\c}$ whenever $J_-+J_+\geq\c$ (resulting with 
the length of $\widetilde{I}$ having the $\Gamma_{\c}(2,1)$ law).
We then remove from the list of active points all those points which  
got encompassed by the interval $I_k$, including the base point
$(w_k,t)$. The links in the graph connected 
to all such points other than $(w_k,t)$, are considered to be surplus edges.

\noindent
\emph{Connections of $I_k$:} 
With $I_k = \{w_k\} \times \widetilde{I}_k$, for each $i \ne w_k$ we
view $\widetilde{I}_k$ as a subset of $\S^{w_k,i}$ and sequentially 
for $i=1,2,\ldots,N$, sample the process of links 
$\calL_{w_k,i}$ for $t' \ge 0$ restricted to $\widetilde{I}_k$. 
We erase all links between $I_k$ and points on already explored intervals, 
and register each link end $(i,t')$ on the neutral space as an active point, 
labeled with the time (order) of its registration.

\noindent
When done examining all the connections from $I_k$, we change its 
status from neutral to that of an explored interval and increase 
$k$ by one, continuing with this procedure till no neutral space remains
(which happens after finitely many steps, since the number of intervals in 
the \abbr{QRG} is finite). To recover the resulting \abbr{QRG} we need only
to keep track of the explored intervals (end-points), and the $\zeta_k$ 
new links that have been formed in each step. 

Now, let $\eta_k=\zeta_k-(\textrm{sur}(k)-\textrm{sur}(k-1))$, where
\textrm{sur}(k) counts all the surplus edges found by the end of each of the first $k$ steps
of exploration.
Then, by definition
\begin{equation}
\label{eq:active:walk}
A_k=\begin{cases}
A_{k-1}+\eta_k -1, \quad &\text{if $A_{k-1}>0$}\\
\eta_k, \quad &\text{if $A_{k-1}=0$\,.} 
\end{cases}
\end{equation}
As mentioned before, the exploration of the first component containing the point $(1,t)$ sampled at the initial stage ends
at $\tau_1 = \min\{k \ge 1: A_k = 0\}$, with its size $|\calC (1,t)|$ being 
$\tau_1$ (the number of explored intervals thus far). A new component whose
size is $\tau_2-\tau_1$ is then explored from step $\tau_1$ 
till the end of step $\tau_2 = \min\{k > \tau_1: A_k =0\}$, and so on.

\subsection{Overcounting}
\label{subsec:overcounting}

Let $m_i^{(T)} \le m_i$ count the intervals in vertex $i$ which 
belong to components of $\calG_N^\beta$ whose sizes exceed $T$
and ${\calC}(i,{\ast})$ denotes the connected component of $\calG^\be_N$ 
containing $\S^i_\c$, \emph{after erasing} all the holes 
punched in $\S^i_\c$ by $\calH_i$. Since the size of the component 
containing interval $I_i^l$ of \eqref{holedecomp} is at most $|\calC (i,\ast)|+m_i-1$, 
it follows that $m_i^{(2T)}=0$ whenever both $m_i \le T$ and $|\calC(i,\ast)| \le T$.
Hence, by Markov's inequality
\begin{align}\label{eqn:simplify-m}
\P\left(|\calC_{\max}| \ge 2T\right) &\le \P\big(\sum_{i=1}^N m_i^{(2T)} \ge 2T \big)\nonumber \\
&\le \frac{1}{2T} \sum_{i=1}^N \E [ m_i^{(2T)} ] \le 
\frac{1}{2T} \sum_{i=1}^{N} \E \Big[ m_i ({\bf 1}_{|\calC(i,\ast)|>T} 
+ {\bf 1}_{m_i > T}) \Big] \,.
\end{align}
Further, $|\calC(i,\ast)|$, $i\in G_N$, are identically distributed random 
variables, each of which is independent of the corresponding variable $m_i$ which in 
turn has the Poisson($\c) \vee 1$ distribution
(so $\E [m_i] = \c + e^{-\c}$). Consequently,
\begin{equation} \label{eqn:simplify-m-1}
\P\left(|\calC_{\max}| \ge 2T\right) \le 
\frac{N (\c+{\rm e}^{-\c})}{2T} \Big[ \P (|\calC(1,\ast)|>T) + \P(m_1 \ge T) \Big] \,.
\end{equation}
Taking $T=(A/2) N^{2/3}=H^2$ we thus establish part (a) of Theorem 
\ref{thm:crit:size}, upon showing that 
for some $c$
 finite and all $H$ and $N$ large enough
\begin{equation}\label{eq:bd-c1-ast}
\P(|\calC(1,\ast)|>H^2) \le \frac{c}{H} \,.
\end{equation}
Since $\calC(1,\ast)$ corresponds to the exploration process starting at 
$I_1=\{1\} \times \S_\c$,
it suffices to consider the value of $\tau_1$ when the corresponding
$\{\eta_k\}_{k \ge 1}$ 
are replaced in \eqref{eq:active:walk} by another, simpler to analyze, collection 
$\{\xi_k\}$ that stochastically dominate them.

For us to be able to estimate tail probabilities of $\tau_1$ using the i.i.d. sequence $\{\xi_k\}$, 
we must define appropriate coupling between $\{\eta_k\}$ and $\{\xi_k\}$. To this end, let us define
\begin{equation}\label{eq:rw-overcount}
S_n = S_1 + \sum_{k=2}^n (\xi_k -1) \,, 
\end{equation}
and consider the following monotone coupling between $A_k$ associated to $\calC(1,\ast)$ 
and $S_k$ up to time $\tau_1$.

Like $A_0$, we begin with setting $S_0=1$. Since $\calC(1,\ast)$ corresponds to the exploration 
process starting at $I_1=\{1\} \times \S_\c$, $\eta_1$ follows Poisson$\left( \frac{(N-1)\c}{N\lambda}\right)$ 
distribution. Let $\xi_1$ be $\eta_1$ together with {\it self-links} of the interval to itself. Therefore,
$\xi_1$ follows Poisson$\left( \frac{\c}{\lambda}\right)$. Since the coupling is up to time $\tau_1$ we only have to consider the 
$A_{k-1} > 0$ case at step $k\geq 2$. As explained above, we choose an active point 
$(w_k,t)$ and sample links included in the counting towards $\zeta_k$. In order to define $\xi_k$ 
recall the procedure of sampling connections and consider the following addition to it. In addition to 
carving out of the maximal neutral interval $I_k$ around $(w_k,t)$, consider also 
the full interval $I_k^S$ around $(w_k,t)$ having the length law $\Gamma_{\c}(2,1)$. In addition to 
the links sampled by $\calL_{w_k,i}$-s $\forall i\neq w_k$ restricted to $\widetilde{I}_k$, run a
unit intensity Poisson process on $I_k^S\cap I_k^c$, and another independent Poisson process of $\left(\frac1N\right)$
intensity on the interval $\widetilde{I}_k$ counting self-links to the same interval. Let 
$\upsilon_k$ count the arrival points of these additional Poisson processes. Define 
$\xi_k=\zeta^S_k+\upsilon_k$ with $\zeta^S_k$ counting all the links created by $\calL_{w_k,i}$-s 
restricted to $\widetilde{I}_k$ without erasing those whose end points fall on already explored 
intervals. Obviously, $\xi_k\geq \eta_k$ for all $k\ge 1$. Note also that for $k\ge 2$, the random variables $\xi_k$ 
are i.i.d. each following Poisson$\left(\frac{J_k}{\lambda}\right)$ conditioned on $J_k$, which are i.i.d. 
$\Gamma_{\c}(2,1)$, and independent of the $\xi_1$.   

Defining $\tau = \min\{n \ge 1: S_n = 0\}$ as the first hitting time of zero by the process $S_k$,
and using the monotone coupling argument, the inequality in \eqref{eq:bd-c1-ast} follows from the bound 
\begin{equation}\label{eq:kemperman-conseq}
\P\left(\tau>H^2\right)\le\frac{c}{H}.
\end{equation}

Having $\{-1,0,1,2,\ldots\}$-valued increments, recall Kemperman's 
formula for such a random walk, stating that for any $\ell \ge 0$ and $n \ge 1$,  
$$
\P(\tau=n+1|S_1=\ell) = \frac{\ell}{n} \P(\left.S_{n+1} = 0\right| S_1=\ell) 
$$ 
(see \cite[Theorem 7, p.165]{GrimStir}). 
Then, we can write
\begin{align}\label{eqn:localclt}
\P\left(\tau=n+1\right) &= \sum_{\ell=0}^n \frac{\ell}{n} \P(S_{n+1}-S_1=-\ell)
\P\left( S_1 = \ell\right) \nonumber \\
& \le \frac{\E(S_1)}{n} \sup_{\ell} \{\P(S_{n+1}-S_1=-\ell)\} \,.
\end{align}
%
Our assumption that $F(\beta,\lambda)=1$, implies 
that $\E \xi_2 = \la^{-1}\E J_2 =1$, so $\{S_n\}_{n\ge 2}$ has zero-mean i.i.d. 
increments of finite exponential tails. Thus, applying the 
local \textsc{clt} for the lattice random walk $S_{n+1}-S_1$
(see \cite[Proposition 2.4.4]{LawlerLimicBook}), we deduce from \eqref{eqn:localclt} that 
\begin{equation}\label{eq:tau-density}
\P(\tau=n+1) \le  c n^{-3/2}\,,
\end{equation}
for some $c$ finite and all $n$, which together with \eqref{eq:kemperman-conseq}, proves that
\begin{equation}\label{eq:cmax-ubd-1}
\P\left( |\calC_{\rm max}| > A N^{2/3}\right) \le c_* A^{-3/2}.
\end{equation}


\subsection{Undercounting}
\label{subsec:free-exploration} 

To bound the lower tail of $|\calC_{\max}|$, we construct a stochastic
lower bound for all component sizes by following a more 
restrictive exploration process, which after forming the  
first active point on each vertex $w \in G_N$, voids all space
on that same vertex beyond the relevant interval around this
active point (thus sequentially producing components with no 
more intervals than does the original exploration process). Specifically,
after the initial stage, at each step $k \ge 1$ the restrictive exploration
considers for $I_k$ only active intervals or completely neutral circles 
(as in parts (a) and (b) of the original exploration process defined in Section \ref{subsec:exploration}), 
till none such are left. It also keeps at most one connection 
from $I_k$ to any, as of yet, never visited 
(in particular, completely neutral) circle $\S^i_\c$, 
ignoring (erasing) all the other links which are being formed 
in step $k$ by the original exploration process.  
Note that this restrictive process has no surplus edges and its
number of active points $A^f_k$, starts at $A^f_0=1$ and follows the
recursion 
\begin{equation}
\label{eq:active:walk-res}
A_k^f=\begin{cases}
A^f_{k-1}+\eta^f_k -1, \quad &\text{if $A^f_{k-1}>0$}\\
\eta^f_k, \quad &\text{if $A^f_{k-1}=0$\,.} 
\end{cases}
\end{equation}
Here, conditioned on $A^f_{k-1}$ and $J_k$, the variables $\eta^f_k$
are independent variables distributed as Bin$(N^f_{k-1},1- e^{-J_k/(\la N)})$
for i.i.d. $\Gamma_\c (2,1)$-distributed collection $\{J_k\}$ and
$N^f_k:=N-k-(A^f_k \vee 1)$. As before, the component sizes are given by
$\tau^f_{r}-\tau^f_{r-1}$, for successive returns to zero 
$\tau^f_r=\min\{k>\tau^f_{r-1} : A_k^f = 0 \}$, starting at $\tau^f_0=0$.

Had we replaced $J_k$ by $\E J_k = \la$, it would have resulted in the exploration 
process for the (effectively) critical Erd\H{o}s-R\'{e}nyi random graph 
$G(N,1-e^{-1/N})$, for which \eqref{eq:cmax-lbd} is well-known, for example, 
see \cite[Theorem 2]{nach-peres1}. As we are not aware of a study of component 
sizes for our inhomogeneous graph, we next adapt the proof of 
\cite[Theorem 2]{nach-peres1} to our context.

First, from the recursion \eqref{eq:active:walk-res} conditioned on the event $\{A_{k-1}^f > 0\}$, then 
$$
(A_k^f)^2-(A_{k-1}^f)^2 = (\eta^f_k-1)^2 + 2(\eta^f_k-1)A^f_{k-1}.
$$
Conditioned on the event $\{0<A^f_{k-1}\le h\}$ for some arbitrary $h>0$, which we shall specify later, we observe that
\begin{equation} \label{m:eqn:free:quadratic}
\E\left[(A^f_k)^2-(A^f_{k-1})^2 \big| A^f_{k-1} \right]
\ge \frac{ (N-h-k)^2}{(\la N)^2} \E[J_k^2] -2 \frac hN\left(h+k \right)+ {\rm O}\left(\frac 1N \right) \,.
\end{equation}
Further, since $\E[J_k^2] > \la^2$, so for $h = c_h N^{1/3}$, for all $k \le T_h$ with $T_h = c_1 N^{2/3}$,
and large enough $N$
\begin{equation}
\label{m:eqn:submart}
\E\left[(A^f_k)^2-(A^f_{k-1})^2 \big| A^f_{k-1} \right] \ge  1- 2 c_h c_1 \,,
\end{equation}
where $c_h$ and $c_1$ are arbitrary positive constants to be chosen later.

The latter bound applies also when $A^f_{k-1}=0$, as then $A^f_k=\eta_k^f$.
Now, taking $b:=1-2 c_h c_1 > 0$, we consider the stopping time 
$$
\tau_h = T_h \wedge \min\{k \ge 0 : A_k^f \ge h\}
$$
noting that by the preceding calculation,
$L_k := (A^f_{k\wedge\tau_h})^2 - b (k\wedge\tau_h)$ is a sub-martingale
starting at $L_0=1$. Further, it is shown in \cite[Proof of Lemma 5]{nach-peres1} 
that for $\xi$ a Bin$(n,p)$ variable, and any $n \le N$, the
distribution of $\xi-r$, conditioned on the event $\{\xi \ge r\}$, is stochastically dominated by Bin$(N,p)$.
Thus, in our setting, given $\tau_h=k \le T_h$, $N^f_{k-1} = n \le N$, 
$A^f_{k-1}=\ell$ and $p_N=1-e^{-J_k/(\la N)}$, with $\{J_k\}$ i.i.d $\Gamma_{\c}(2,1)$, we have that conditioned on the event
$\{A^f_k \ge h\}$, the distribution of $(A^f_k-h)$ is stochastically 
dominated by Bin$(N,p_N)$ conditioned on the same $p_N$. 
Averaging over all possible $n,k,\ell,p_N$
values, we deduce that conditioned on 
$\left\{ A^f_{\tau_h} \ge h\right\}$ the over shoot $\left(A^f_{\tau_h} - h \right)$ 
is stochastically dominated by Bin$(N,p_N)$ conditioned on 
$p_N=1-e^{-J/(\la N)}$ with $J \sim \Gamma_{\c}(2,1)$. Consequently,
$$
\E[(A^f_{\tau_h})^2] \le h^2 + 2 h \E[N p_N] + \E[ (N p_N)^2 ] + \E[ N p_N(1-p_N)] 
\le (h+2)^2 \,,
$$
for $h = c_h N^{1/3}$ and all $N$ large enough. With $\tau_h \le T_h$, 
upon applying the optional stopping theorem for the $L^2$-bounded 
sub-martingale $\{L_k\}$, 
we find that 
$\E[(A^f_{\tau_h})^2] \ge 1 + b \E(\tau_h)$. Thus, by Markov's inequality,
\begin{equation} \label{m:eqn:tau-h}
\P(\tau_h = T_h) \le \frac{\E[\tau_h]}{T_h} \le \frac{(h+2)^2}{b T_h} 
= \frac{c_h^2}{b c_1} (1+o(1))\,.
\end{equation}
Fixing $T_0 = \delta N^{2/3}$ note that $|\calC_{\rm \max}|$ exceeds the
value of the stopping time
$$
\tau_0 = T_0 \wedge \min\{s \ge 0: A^f_{\tau_h+s} =0\}\,,
$$ 
so the stated bound \eqref{eq:cmax-lbd} follows once we show that   
for any $b' > \la^{-2} \E [J^2] + 2 c_h (c_1+\delta)$,
\begin{equation}\label{m:eq:tau-o-bd}
\P\left(\left.\tau_0 < T_0\,\right\vert\, A^f_{\tau_h} \ge h\right) 
\le \frac{b' \delta}{c^2_h} \,.
\end{equation}
Indeed, we then 
choose $c_h = \frac{1}{2} \delta^{1/5}$ and $c_1=1/(4 c_h)$, so $b=1/2$
and $b'=\beta^2+1$ works whenever $\delta \le 1/2$. 

To derive \eqref{m:eq:tau-o-bd} consider the uniformly 
bounded, non-negative process $M_k = \max\{h-A^f_{\tau_h+k},0\}$.
If $0 < M_{k-1} < h$ then by \eqref{eq:active:walk-res}
$$
M^2_k - M^2_{k-1} \le \left(\eta^f_{\tau_h+k} -1 \right)^2 + 2\left( 1-\eta^f_{\tau_h+k}\right) M_{k-1} \,.
$$
The same inequality applies when $M_{k-1}=0$ (i.e. $A^f_{\tau_h+k-1} \ge h$, 
so $M_k \le \max\{1-\eta^f_{\tau_h+k},0\}$). By definition, 
$A^f_{\tau_h+k-1} \le h$ whenever $M_{k-1} \ne 0$, hence
for $N$ large enough and all $k \le T_0$, $\tau_h \le T_h$ we 
find, as in the derivation of \eqref{m:eqn:submart} that  
$$
\E\Big[ M^2_k - M^2_{k-1}\,\vert\, M_{k-1}<h \Big] \le 
\la^{-2} \E[J^2] + 2 \frac{h}{N} (T_h + T_0 + h) + o(1) \le b' \,.
$$ 
Thus, conditioned on the event $A^f_{\tau_h} \ge h$ the process 
$\{M^2_{k \wedge \tau_0} - b' (k\wedge \tau_0)\}$ 
is a super-martingale which starts at zero. 
Noting that $\{ \tau_0<T_0 \} \subseteq \{M_{\tau_0}=h\}$, upon 
applying the optional stopping theorem for this process at $\tau_0$,
we conclude that
$$
\P\left(\left.\tau_0 < T_0\,\right\vert\, A^f_{\tau_h} \ge h\right) 
\le h^{-2} \E\big[  M^2_{\tau_0} \,\vert\, A^f_{\tau_h} \ge h\big] 
\le b' h^{-2} \E [\tau_0\,\vert\, A^f_{\tau_h} \ge h ] \le b' \frac{T_0}{h^2} \,,
$$
as stated. \qed


\begin{proof}[Proof of Theorem \ref{thm:crit:size}] Observe that $\calP(\calC_{\max})$ is stochastically dominated by 
$\sum_{i=1}^{|\calC_{\max}|}\xi_i$,
where $\{\xi_i\}$ are i.i.d. cut-Gamma random variables. Repeating the arguments set forth in Section \ref{subsec:overcounting}
we can easily conclude the upper bound as in equation \eqref{eq:cmax-ubd}.

Similarly, we shall propose a process using the same undercounting algorithm such that it is stochastically dominated by
$\calP(\calC_{\max})$. In particular, instead of accounting only the number of once visited vertices, let us associate a cut 
gamma $\Gamma_{\c}(2,1)$ random variable with every such vertex visited. 
Then consider the sum of all such random variables, which clearly is dominated by
$\calP(\calC_{\max})$. Thereafter, again following the same steps as in Section \ref{subsec:free-exploration} we 
conclude the required result of Theorem \ref{thm:crit:size}.
\end{proof}

\section{Proof of Theorem \ref{thm:joint:cgs}}
\label{sec:joint:cgs}

\subsection{Exploration process and Brownian excursions}
\label{subsec:brownian:excursions}

Recall the 
length of a sampled interval of the \abbr{qrg} being 
$J=\min(J_-+J_+,\c)$ for i.i.d. standard Exponential 
variables $J_-,J_+$ (a distribution we denote by $\Gamma_\c (2,1)$).
With 
$\E(J)=F(\c)$ of \eqref{eq:critical}, the critical curve has the explicit expression  
$\lambda=F(\c)$. Further, the critical window around some $\lambda_\star=F(\c_\star)$
for $\c_\star>0$, corresponds to fixing $a \in \real$ and considering
\begin{equation}\label{eq:window-curve}
\c_N \to \c_\star \,, \qquad F(\c_N) = \la_N \big(1+a N^{-1/3}\big) \,.
\end{equation}

Let $\widetilde{N} = \sum_{i=1}^N m_i$ denote the total number of steps in 
the exploration process of Section \ref{subsec:exploration},
and $(Y_N^{a,\cn}(k), k \le \widetilde{N})$ be the 
breadth-first walk associated with the \abbr{qrg} on $\calG_{\c_N}^N$, where 
$(\c_N,\lambda_N)$ satisfy \eqref{eq:window-curve}. That is,
\begin{equation}\label{eq:YN-def}
Y_N^{a,\cn}(k)=Y_N^{a,\cn}(k-1)+\eta_k-1, \;\; \qquad
Y_N^{a,\cn}(0)=1,
\end{equation}
for $\eta_k$ of recursion \eqref{eq:active:walk}. Thus, $Y_N^{a,\cn}(k)$ 
(which may well become negative as $k$ grows), counts the number of active points 
at the end of step $k$, minus the number of explored components before step $k$. 

As in \cite{aldous}, observe that lengths of excursions of the process $Y^{a,\c_N,+}_N(k) = 
Y^{a,\c_N}_N(k) - \min_{l\le k}Y^{a,\c_N}_N(l)$ above zero correspond to size of the connected 
component containing the vertex where the process $Y^{a,\c_N}_N(k)$ started.

Setting 
\begin{equation}\label{eq:sig-star-def}
\sigma_\star^2 = \sigma^2(\c_\star)
 = \frac{1}{\lambda_*^2} \E\Big[ 
\big( J_-+J_+) \wedge \c_\star\big)^2 \Big]\,,
\end{equation}
our goal in this subsection is to prove the following proposition.
\begin{proposition}
\label{thm:scaling:limit}
For $(\c_N,\lambda_N)$ that satisfy \eqref{eq:window-curve},
as $N \to \infty$, the processes
\begin{equation}
\label{eq:scaling}
\bar{Y}_N^{a,\cn}(s) =  
N^{-1/3} Y_N^{a,\cn}(\lfloor N^{2/3}s\rfloor \wedge \widetilde{N})\,,
\end{equation}
converge in law to $W^{a,\c_\star}$ of \eqref{eqn:def:W-a-beta}  
(on the space $D\left( [0,\infty)\right)$ 
equipped with the topology of uniform convergence on compacts). 
\end{proposition}  

Recall that 
$\zeta_k$ 
links are generated at step $k$ of the exploration process and let
\begin{equation}\label{eq:ZN-def}
Z^{a,\cn}_N(l) = Z^{a,\cn}_N(l-1) + \zeta_l -1\,, \qquad Z^{a,\cn}_N(0) =1,
\end{equation}
be the corresponding breadth-first walk. Since we sample intervals only 
when they are to be explored, the walk $Z_N$ does not distinguish 
between active points that end as intervals of the \abbr{QRG} and 
those that are later found to be on surplus edges. Nevertheless, our 
next proposition controls the number of active points, which 
as seen in Remark \ref{rem:sur-bd}, yields having at most ${\rm O}_\P(1)$
surplus edges till step $s N^{2/3}$.
\begin{proposition}
\label{prop:nu-1-0}
Fixing $a\in\real$, recall the count $A_k$ of active points at the end of step $k$ 
of the exploration for $(\c_N,\la_N)$ satisfying \eqref{eq:window-curve}. 
Then, for some $K=K(a,s)<\infty$ and all $L,N \ge L_0(a,s)$, 
\begin{equation}
\label{eqn:claim:order}
\P\left(N^{-1/3}\, \max_{k \le s N^{2/3}} \{ A_k \} > L\right) \le K L^{-2}\,.
\end{equation}
\end{proposition}

The proof of above proposition involves elaborate, but crude, bounds
on functionals of $Z_N$, and thus we defer it to the appendix.

\begin{remark}\label{rem:sur-bd}
Recall $\textrm{sur}(l)$ counts the surplus links detected
in part (a) of the exploration process during its first $l$ steps. 
The order of exploring active points is such that 
the first active point formed on any given vertex never 
contributes to $\textrm{sur}(l)$. Further, $\textrm{sur}(l)$
is bounded above by the aggregate count $\Delta_{\rm e}(l)$ 
of active points on vertex $w_k$ at step $k \le l$,
beyond the explored point $(w_k,t)$.
Conditional on the state of the process
at the start of its $k$-th step, the number of active points registered
during that step that may contribute to $\Delta_e(l)$, 
is stochastically dominated by a 
Poisson$(A_{k-1} \cn (\lambda_N  \, N)^{-1})$ random variable. In particular,
for some $\kappa$ finite and all $N$, $l$,
\begin{equation}
\E [ \textrm{sur} (l) ] \le 
\E [ \Delta_{\rm e} (l) ] \le
\frac{\kappa}{N} \E \big[\,\sum_{k=0}^{l-1} A_k\,\big]
\le \frac{\kappa l}{N} \E \big[\,\max_{k \le l} \{A_k\} \,\big] \,.
\end{equation} 
From Proposition \ref{prop:nu-1-0} we have 
that $\E[ \Delta_{\rm e}(s N^{2/3})]$ and hence
$\E [ \textrm{sur} (s N^{2/3})]$, is uniformly bounded in $N$. 
\end{remark}

The control on number of vertices which the exploration process 
visits at least twice by the end of the $k$-th step 
(for $k=s N^{2/3}$), is crucial for the success of our analysis. 
To this end, we define hereafter \emph{the
number of visits to vertex $v \in G_N$ by the end of the $k$-th 
exploration step}, as the total number 
of active points formed on $\S_{\cn}^v$ by that time, i.e., 
we count past active points which were removed from the list 
and also those which are currently active by the end of time $k$.
\begin{proposition}
\label{prop:nu-m}
For the exploration of the \abbr{qrg} at parameters 
satisfying \eqref{eq:window-curve}, 
let $\nu^l_{\ge m}$ count the total number of visits 
by the end of its $l$-th step, to sites 
(circles), having at least $m$ such visits each. 
Then, for some finite $\kappa$, all positive $A,s$ and $N$ 
large enough, we have that for $m=1,2,3$ and any $l \in [1,s N^{2/3}]$, 
\begin{equation}\label{eq:num-event-bd}
\P ( \nu_{\ge m}^l \, \ge \, A l^m N^{1-m}) \le \frac{\kappa^m}{A} \,. 
\end{equation}
\end{proposition}
 
\begin{proof} Let $\calF_t$ denote the filtration generated 
by the state of the exploration process of Subsection 
\ref{subsec:exploration},
namely, the neutral sub-space and collection of active points,
with $\calF_k$ for integer $k \ge 0$, 
denoting the state at the end of step $k$. Further, let 
$\tau(v)$ the stage in which 
it first visits $v \in G_N$ (so $\calF_{\tau(v)}$ 
records the state of the exploration process 
immediately after selecting its first active point on $\S_\c^v$). 
For any $l \ge 1$ and $v$ let $L_v(l)$  
count the links whose end points are on the neutral part of $v$\footnote{
We note here that
$L_v(l)$ is different from the total number of visits to the vertex $v$ by time $l$
because according to our description, the point where we start the exploration process,
and the points where we restart our exploration process after $A_k = 0$ are indeed
counted as visits but these are not identified as end points of links.}
till the end of the $l$-th step of that process, with $L'_v(l)$
counting only such links made after $\tau(v)$ 
(setting $L_v'(l)=0$ in case $\tau (v) \ge l+1$). Only 
one interval is explored in each step, hence
\begin{equation}\label{eq:nu1-exp}
\E [ \nu^l_{\ge 1} ] \le l + \sum_{v=1}^N \E [L_{v}(l)] \,.
\end{equation}
Note that $L_{v}(l)$ increases when considering the overcounting process,
so the $l$ explored intervals are complete circles
(of length $\cn$), other than the circle at $v$, which remains completely neutral, even after links to it are formed
and intervals are sampled around the links.
Thus, $L_{v}(l)$ is stochastically dominated by a Poisson random
variable with parameter 
\begin{equation}\label{eq:dom-Poi}
\frac{\cn\, l}{\lambda_N\, N} \le \frac{\kappa' l}{N} \,,
\end{equation}
for some $\kappa'$ finite and all $N$ large enough. In particular,
by \eqref{eq:nu1-exp} $\E[\nu^l_{\ge 1}] \le (\kappa'+1) l$, 
which in combination with Markov's inequality 
establishes \eqref{eq:num-event-bd} for $m=1$ and any $\kappa \ge \kappa'+1$.

Next recall that for $k \le N$, either $A_k>0$ so part (a) 
of the exploration process  
applies at the $k$-th step, or else part (b) applies for it
(since at most $k-1<N$ vertices have been explored before).
Consequently, assuming hereafter that $s N^{2/3} \le N$, part (c)
of the exploration does not occur throughout its first $l$ steps. 
Further, all active points chosen in part (b) or the initial stage 
of the process result with a \emph{first} visit of new vertex. 
Hence we have in analogy with \eqref{eq:nu1-exp} that 
\begin{equation}\label{eq:num-basic-bd}
\E[ \nu^l_{\ge 2} \,] \le 2 \sum_{v=1}^N \E [ L'_{v}(l) ]  \,.
\end{equation}
As argued before, conditional on $\calF_{\tau(v)}$ the value of $L'_{v}(l)$ 
increases if from time $\tau(v)$ onward we modify the process to have
all explored intervals be complete circles (of length $\cn$), on 
vertices other than $v$, while keeping the circle at $v$ completely neutral. 
That is, conditionally on $\calF_{\tau(v)}$
the variable $L_v'(l)$ is stochastically dominated by a Poisson 
variable of parameter $\kappa' l N^{-1}$ times the indicator 
on the event $\{\tau(v) < l+1 \}$. Hence, for any $l \ge 1$ and $v \in G_N$,
\begin{align}
\E \big[L'_{v}(l) \, \big | \, \calF_{\tau(v)} \big]
\le \Big(\frac{\kappa' l}{N}\Big) \, {\bf 1}_{\{\tau(v) < l+1\}} \,. 
\label{eq:num-iterate}
\end{align}
Summing over $v$ the expected value of \eqref{eq:num-iterate}, 
we deduce from \eqref{eq:num-basic-bd} that 
\begin{equation}\label{eq:num-exp}
\E[ \nu^l_{\ge 2} \,] \le 2 \Big(\frac{\kappa' l}{N}\Big) \, 
\E[\nu^l_{\ge 1}] \le 2 (\kappa'+1) \kappa' l^2 N^{-1} 
\end{equation}
from which we recover 
\eqref{eq:num-event-bd} for $m=2$ and $\kappa = 2 \kappa'+1 $ 
(by Markov's inequality).

Finally, repeating this argument, now with $\tau(v)$ 
the time at which the second active point on $v \in G_N$ 
is selected,
we deduce that 
\begin{equation*}
\E[ \nu^l_{\ge 3} \,] \le 3 \Big(\frac{\kappa' l}{N}\Big) \, \E[\nu^l_{\ge 2} \,],
\end{equation*}
which upon suitably increasing the value of $\kappa$, yields
\eqref{eq:num-event-bd} for $m=3$.
\end{proof}


\begin{proof}[Proof of Proposition \ref{thm:scaling:limit}] Equipping 
$D([0,\infty))$ with the topology
of uniform convergence on compacts, let
\begin{equation}\label{eq:scaling-Z}
\bar{Z}_N^{a,\cn}(s)=N^{-1/3} Z_N^{a,\cn}(\lfloor N^{2/3}s\rfloor \wedge \widetilde{N})\,,
\end{equation}
for the breadth first walk $Z_N^{a,\cn}(\cdot)$ of \eqref{eq:ZN-def}. Recall that 
$$
\bar{Z}_N^{a,\cn}(s) - \bar{Y}_N^{a,\cn}(s) = N^{-1/3} \textrm{sur}(
\lfloor N^{2/3}s\rfloor \wedge \widetilde{N}) 
$$
is non-decreasing in $s$ and so by Remark \ref{rem:sur-bd}, as $N \to \infty$,
\begin{equation}\label{eqn:first-step}
\sup_{s\le s_0} \left|\bar{Z}_N^{a,\cn}(s) - \bar{Y}_N^{a,\cn}(s) \right| \rightarrow 0 \,\,\,\ \ \ \text{ in }\P\,.
\end{equation}
It thus suffices to prove that $\bar{Z}_N^{a,\cn}$
converges in law to the desired limit $W^{a,\c_\star}$.
To this end, by Doob's decomposition with respect to the canonical 
filtration $\calF_k$ associated with the exploration process, 
we get that
\begin{align}\label{eq:dec12}
Z_N^{a,\cn} = M_N^{a,\cn}+ B_N^{a,\cn}\,,\qquad
\big( M_N^{a,\cn} \big)^2 = Q_N^{a,\cn}+D_N^{a,\cn} \,,
\end{align}
with martingales $M_N^{a,\cn}$, $Q_N^{a,\cn}$ (null at $k=0$),
and predictable processes $B_N^{a,\cn}$ and $D_N^{a,\cn}$.
Adopting the notation $\bar{M}_N$, $\bar{B}_N$, 
in accordance with \eqref{eq:scaling-Z}, 
and $\bar{Q}_N$, $\bar{D}_N$ similarly scaled by extra factor 
$N^{-1/3}$ in accordance to the \abbr{rhs} of \eqref{eq:dec12},
we show in Lemmas \ref{lemma:D-N}
and \ref{lem:M_bar}, respectively, that for 
$\sigma_\star$ of \eqref{eq:sig-star-def} and
any finite $s_0$, as $N\rightarrow \infty$,
\begin{align}
\label{eq:D_bar}
& \bar{D}_N^{a,\cn}(s_0)\stackrel{\P}{\to} \sigma_\star^2 s_0\,,\\
\label{eq:M_bar}
& \mathbb{E} \Big[\sup_{s\leq s_0}\big|\bar{M}_N^{a,\cn}(s)-
\bar{M}_N^{a,\cn}(s^-) \big|^2 \Big] \rightarrow 0 \,.
\end{align}
Combining \eqref{eq:D_bar} and \eqref{eq:M_bar},
it then follows from \cite[Theorem 7.1.4]{EthierKurtz} that 
the martingales $\{\bar{M}_N^{a,\cn}\}$ 
converge weakly in $D([0,\infty))$ to $\sigma_\star W$ for a
standard Brownian motion $W$.
Further, we show in Proposition \ref{lemma:B_bar} that 
\begin{equation}
\label{eq:B_bar}
 \sup_{s\leq s_0}\left|\bar{B}_N^{a,\cn}(s)-\rho^{a,\c_\star} (s)\right|
\rightarrow_{\P} 0\,. 
\end{equation}
That is, the sequence of predictable 
processes $\bar{B}_N^{a,\cn}$ converges in probability 
in $D([0,\infty))$ to the non-random $\rho^{a,\c_\star}$
of \eqref{eqn:def:rho-a-beta}, hence $\bar{Z}_N^{a,\cn}$ 
converges in law to $W^{a,\c_\star}$.
\end{proof}

Proceeding with the proof of \eqref{eq:D_bar}--\eqref{eq:B_bar},  
we often drop the indices $(a,\cn)$ from $Z_N^{a,\cn}$ and related 
random variables. We start by establishing \eqref{eq:D_bar}. That is,
\begin{lemma}
\label{lemma:D-N} For $(\cn,\la_N)$ that satisfy \eqref{eq:window-curve}, 
$\sigma_\star$ of \eqref{eq:sig-star-def}, any $\delta>0$ and $s_0$ finite,  
\begin{equation} \label{eqn:lemma:D-N}
\lim_{N \to \infty} \,
\sup_{l \le N^{2/3} s_0} \, \P ( |D_N^{a,\cn}(l) - \sigma_\star^2 l| 
\ge 2 \delta N^{2/3}) = 0 \,. 
\end{equation}
\end{lemma}
\begin{proof} Recall that $D_N(k) - D_N(k-1) = {\rm var} (\zeta_k|\calF_{k-1})$
(starting at $D_N(0)=0$). 
Hereafter, $\calF_{k^-}$ denotes $\calF_{k-1}$ augmented by 
the active point $(w_k,X_k)$ of step $k$ and the 
interval $I_k=\{w_k\} \times \widetilde{I}_k$ around it. The law of 
$\zeta_k$ given $\calF_{k^-}$
is Poisson$(\varphi_k)$, with $\varphi_k$ denoting
the aggregate over circles other than $\S_{\cn}^{w_k}$, of the length 
of their neutral space restricted to $\widetilde{I}_k$ and divided by $\la_N N$. 
As such, both its conditional mean and conditional variance
are given by $\varphi_k$, hence by the variance conditioning decomposition
(at $\calF_{k^-}$), 
\begin{equation}\label{eq:zeta-to-phi}
{\rm var} (\zeta_k|\calF_{k-1}) = \E [ \varphi_k | \calF_{k-1} ] + 
{\rm var} [\varphi_k | \calF_{k-1} ] \,.
\end{equation}
Further, at the start of the $k$-th step there are at least $N-k$ 
completely neutral circles beyond the vertex on which the $k$-th
explored interval lies. Hence, for $J_k=|\widetilde{I}_k|$ and 
any $k \le s_0 N^{2/3}$,
$$
0 \le \frac{J_k}{\la_N} - \varphi_k \le 
\frac{J_k k}{\la_N N} \le \kappa' N^{-1/3}  
$$
(for $\kappa' = s_0 \sup_N \{\cn/\la_N\}$ finite). It thus suffices to 
prove \eqref{eqn:lemma:D-N} 
for $\widehat{D}_N(l)$ instead of $D_N(l)$, where 
$$
\widehat{D}_N(k)-\widehat{D}_N(k-1) = 
\E [ J_k/\la_N \,|\, \calF_{k-1}] + {\rm var} [J_k/\la_N \, | \, \calF_{k-1} ] =: \Delta_k \,.
$$
The non-negative $\Delta_k$ are uniformly bounded by 
$\overline{\Delta} := \sup_N \{ (\cn/\la_N)+(\cn/\la_N)^2\}$.
Moreover, whenever $\calF_{k-1}$ dictates that the $k$-th step 
explores the first interval on a given vertex $w \in G_N$, it yields
a conditionally independent $J_k$ that follow the $\Gamma_{\cn}(2,1)$ 
distribution. We consequently have that $\Delta_k=F_2(\c_N,\la_N)$
in any such step, where 
\[
F_2(\c,\la) := \E[(J/\la)]+{\rm var}(J/\la) \,.
\]
This applies to all but at most $\nu^l_{\ge 2}$ of the first $l$ steps, 
hence
$$
|\widehat{D}_N(l) - F_2(\beta_N,\la_N) l| \le \overline{\Delta}\, \nu^l_{\ge 2} \,. 
$$
From Proposition \ref{prop:nu-m} we know that 
$\P(\nu^l_{\ge 2} \ge \delta N^{2/3}) \to 0$ for $l=s_0 N^{2/3}$ and
we thus get \eqref{eqn:lemma:D-N} for $\widehat{D}_N(l)$ upon 
noting that $
F_2(\c_\star,\la_\star) = 
\sigma^2_\star$ of \eqref{eq:sig-star-def}.
\end{proof}

We next establish \eqref{eq:M_bar}, thereby
moving closer to completing the proof of
Proposition \ref{thm:scaling:limit}.
\begin{lemma}
\label{lem:M_bar} 
For $(\cn,\la_N)$ that satisfy \eqref{eq:window-curve}
and any $s_0$ finite, 
\begin{equation}\label{eq:M-N-lem}
\lim_{N \to \infty}
N^{-2/3} \mathbb{E} \big[
\max_{1 \le l\leq s_0 N^{2/3}}\, |M^{a,\cn}_N(l)-M^{a,\cn}_N(l-1)|^2
\big] = 0.
\end{equation}
\end{lemma}
\begin{proof} Recall that $M_N(\cdot)$ is the martingale part 
of $Z_N(\cdot)$. Hence, from \eqref{eq:ZN-def} 
$$
M_N(l) - M_N(l-1) = \zeta_l  - \E(\zeta_l | \calF_{l-1})
$$
and \eqref{eq:M-N-lem} amounts to showing that 
$$
N^{-2/3}\E \big[ \max_{l\le s_0 N^{2/3}}
\left( \zeta_l  - \E(\zeta_l | \calF_{l-1}) \right)^2 \big] 
\to 0.
$$
Clearly, 
$\left( \zeta_l  - \E(\zeta_l | \calF_{l-1}) \right)^2 \le 2 \zeta^2_l + 
2 \E(\zeta_l | \calF_{l-1})^2$. Further,  in Section \ref{subsec:overcounting} we saw that conditionally on $\calF_{l-1}$ 
the variable $\zeta_l$ is stochastically dominated by the independent
$\xi_l \ge 0$ whose mean $F(\be_N,\la_N)$ is uniformly bounded (in $N$).
Hence, $\sup_l \E(\zeta_l | \calF_{l-1})^2 \le \sup_N F(\be_N,\la_N)^2$ is finite
and it suffices to show that for i.i.d. $(\xi_l)$,
\begin{equation}\label{eq:bd-zeta-max}
N^{-2/3} \E[ \max_{l \le s_0 N^{2/3}} \zeta^2_l ] \le 
N^{-2/3} \E[ \max_{l \le s_0 N^{2/3}} \xi_l^2 ] \to 0 \,.
\end{equation}
Finally, recall \cite[equation (6')]{HartDav} that the expected maximum 
of $n$ i.i.d. variables of zero-mean and unit variance is 
at most $(n-1)/\sqrt{2n-1}$. Consequently, the expectation on the right 
side of \eqref{eq:bd-zeta-max}
grows at most at rate ${\rm O}(N^{1/3})$, which proves 
\eqref{eq:M-N-lem} (and thereby \eqref{eq:M_bar} as well). 
\end{proof}

\medskip
For the remainder of Section \ref{subsec:brownian:excursions} 
we complete the proof of 
Proposition \ref{thm:scaling:limit} by establishing \eqref{eq:B_bar}. 
Indeed, upon rearranging the expression \eqref{eqn:def:rho-a-beta} for
$\rho^{a,\c_\star}$, this is precisely the 
statement of our next proposition. 
\begin{proposition}\label{lemma:B_bar}
For $(\cn,\la_N)$ that satisfy \eqref{eq:window-curve}
and any $s_0$ finite, as $N \to \infty$, 
\begin{equation}\label{eq:B-lem-bar}
\sup_{s\le s_0} \Big|\bar{B}^{a,\cn}_N(s) - a\,s  
+ \frac{s^2}{2}\Big(1-\frac{F(\c_\star)}{\c_\star} -\gamma(\c_\star) \Big)
+ \frac{s^2}{2}\frac{F(\c_\star)}{\c_\star} \Big|\stackrel{\P}{\to} 0\,.
\end{equation}
\end{proposition}

As the starting point for Proposition \ref{lemma:B_bar}, we
provide the geometric quantities behind the coefficients of $s^2$ in 
\eqref{eq:B-lem-bar}.
\begin{lemma}\label{lem:geo}
For any given interval ${\cal I} \subseteq \S_{\c}$, let
\begin{equation}\label{eqn:H-dfn}
H({\cal I}) := \E [\, |{\cal I}' \cap {\cal I}| \,] \,, \qquad
H_0 ({\cal I}) := \E [\, |{\cal I}' \cap_0 {\cal I}| \,] \,, 
\end{equation}
denote the expectation over an independent interval 
${\cal I}' \subseteq \S_{\c}$ of length law 
$\Gamma_{\c}(2,1)$, built around $0$, where 
\begin{equation}\label{eqn:cap0-dfn}
{\cal I}' \cap_0 {\cal I} := 
\begin{cases} 
\{ x \in {\cal I}' :  
x \not\leftrightarrow 
0 \hbox{ 
within } {\cal I}' \setminus {\cal I} \} \,, \quad & 0 \notin {\cal I} \,, \\
\qquad \emptyset\,,  \quad & 0 \in {\cal I} \,.  
\end{cases}
\end{equation}
For uniform $U \in \S_{\c}$, independent of ${\cal I}$, let
\begin{equation}\label{eq:UH}
\mathsf{UH}(|{\cal I}|)    := \E [ \, H({\cal I}-U) \,| {\cal I}] \,, \qquad
\mathsf{UH}_0(|{\cal I}|)  := \E[ H_0({\cal I}-U) \, | \, {\cal I}] \,.
\end{equation}
Then, for $J$ of $\Gamma_{\c}(2,1)$ law 
\begin{equation}\label{eq:UHav}
\E [\,\mathsf{UH} (J)] = \frac{F(\c)^2}{\c} \,, \qquad
\E [\,\mathsf{UH}_0 (J)] = F(\c) \big(1 - \frac{F(\c)}{\c} - \gamma(\c) \big) \,.
\end{equation}
\end{lemma}
\begin{proof} Recall that for two arcs ${\cal I}$ and ${\cal I}'$ in $\S_{\c}$
of uniformly chosen relative shift $U$, the expected length of 
${\cal I}' \cap ({\cal I}-U)$ is the product of arc lengths divided by $\c$.
In particular, $\mathsf{UH} (J)= \frac{F(\c)}{\c} J$, from which 
the \abbr{lhs} of \eqref{eq:UHav} follows. Similarly, per arc 
${\cal I}'=[-J_-,J_+]$ of length $J$ around $0$ in $\S_{\c}$ 
and $x= \c - |{\cal I}|$, the 
expectation of $\c |{\cal I}' \cap_0 ({\cal I} - U)|$ 
over the uniform shift $U$, is $x (J-x)$ for $J = \c$, 
while for $J<\c$ it is 
\[
\int_{(J_- - x)_+}^{J_-} u \, du + \int_{(J_+ -x)_+}^{J_+} u \, d u 
= x (J - x) + \frac{1}{2} \Big[ 
\big((x-J_+)_+\big)^2 + \big((x-J_-)_+\big)^2\Big] \,,
\]
by elementary geometric considerations.
Computing the expectation of this expression for 
$J_-$ and $J_+$ independent Exp(1) variables, yields
\begin{equation}\label{eq:H0-U-val}
\mathsf{UH}_0 (\c-x) = \frac{1}{\c} \Big[
x F(\c) - \widehat{F}(x) \Big] \,,
\end{equation}
for $\widehat{F}(\cdot)$ as given in \eqref{eqn:gamma-dfn}. Finally, note that 
$\gamma(\c)$ of \eqref{eqn:gamma-dfn} was set so the \abbr{rhs} of \eqref{eq:UHav}
be the expectation of \eqref{eq:H0-U-val} when $(\c-x) \sim 
\Gamma_{\c}(2,1)$.
\end{proof}
Proceeding with the proof of Proposition \ref{lemma:B_bar}, we next 
express $B_N(l)$ as the sum of the terms 
\eqref{eq:bn-expres1}-\eqref{eq:bn-expres3}, which for $l=s N^{2/3}$ 
upon further scaling by $N^{-1/3}$ converge to the three 
limit expressions 
in \eqref{eq:B-lem-bar}, respectively. 
\begin{lemma}\label{lemma:B-N}
Let $X_k$ denote 
the position on $\S^{w_k}_{\cn}$ around which 
$\widetilde{I}_k$
is carved (and if there are no active points 
by the end of exploration step $(k-1)$, then
any function of $X_k$ is replaced by its 
expectation over a uniform $U \in \S^{w_k}_{\cn}$).
Similarly, for $1 \le k < l$ let $\{X_{(k,i)}, 1 \le i \le r(k,l] \}$ 
be the collection of links formed on $\S_{\cn}^{w_{k}}$ during steps $(k,l]$
(some may be in $\widetilde{I}_k$), with the convention that
$r(k,l]=0$ if vertex $w_k$ has been 
explored before. Then, for
$(\cn,\la_N)$ that satisfy \eqref{eq:window-curve}, any
fixed $s_0$ and all $N$ large enough, we have 
uniformly over $l \in [1,s_0 N^{2/3}]$, 
\begin{align}\label{eq:bn-expres1}
B_N(l)=
{\rm O}_{\P}(1)
 & + l \Big(\frac{F(\cn)}{\lambda_N} - 1 \Big)  \\
&- 
\frac{1}{\la_N} \;\;\; \sum_{k=1}^{l} \sum_{i=1}^{r(k,l]} 
 H_0 (\widetilde{I}_{k}-X_{(k,i)})
\label{eq:bn-expres2} \\
&-
\frac1{\lambda_N\,N}\sum_{1 \le k < k' \le l} 
 H(\widetilde{I}_{k}-X_{k'})\, 
 \,.
\label{eq:bn-expres3}
\end{align}
\end{lemma}
\begin{proof}
Recall from proof of Lemma \ref{lemma:D-N} that conditional on ${\calF}_{k'^{-}}$
the number of links $\zeta_{k'}$ formed during the $k'$-th
exploration step is Poisson distributed, whose parameter 
\[
\varphi_{k'} :=\frac{1}{\la_N} |\widetilde{I}_{k'}| - 
\frac{1}{\la_N N} \sum_{k=1}^{k'}  |\widetilde{I}_{k'} \cap \widetilde{I}_{k}| \,,
\]
is bounded by the non-random $\cn/\la_N$. 
Since $B_N(l)$ is the predictable process 
in Doob's decomposition of $Z_N(l)$ we get from \eqref{eq:ZN-def} that 
\begin{align}\label{eqn:B-N:simplify}
B_N(l) &+ (l-1) =  \sum_{k'=1}^{l} \E(\zeta_{k'} | \calF_{k'-1}) =
\sum_{k'=1}^{l} \E (\varphi_{k'} |\calF_{k'-1}) = \frac{1}{\la_N} \Big( l F(\cn) 
\nonumber \\
&- \sum_{k'=1}^{l} \big(F(\cn) - \E (|\widetilde{I}_{k'}| | \calF_{k'-1}) \big)
- \frac{1}{N} \sum_{1 \le k \le k' \le l} 
\E ( |\widetilde{I}_{k'} \cap \widetilde{I}_{k}| | \calF_{k'-1} ) \, \Big) \,.
\end{align}
Hereafter, for $m \ge 1$ let $E^k_m$ and $E^k_{\ge m}$, 
denote the event that upon exploring the interval 
$\widetilde{I}_k$ during the $k$-th step, the site 
$w_k$ on which it lies has been visited precisely 
$m$ times, or at least $m$ times, respectively.
Equipped with these notations we next show that the difference between 
the left sum in \eqref{eqn:B-N:simplify} and the sum in 
\eqref{eq:bn-expres2} is merely part of the 
${\rm O}_{\P}(1)$ in \eqref{eq:bn-expres1}. Specifically, all 
terms are at most $F(\cn)$, which in turn is uniformly bounded 
over $\cn \to \c_\star$, so it suffices to show that there are 
only ${\rm O}_{\P}(1)$ differing terms between these sums.
To this end, since $H_0(\widetilde{I}_k - X_{(k,i)})=0$
whenever $X_{(k,i)} \in \widetilde{I}_k$ (see \eqref{eqn:cap0-dfn}), 
such superfluous links never contribute to \eqref{eq:bn-expres2}. Similarly,
there is no contribution to the left sum in \eqref{eqn:B-N:simplify}
when $E^{k'}_1$ occurs (i.e. when $\widetilde{I}_{k'}$ is the first explored interval on its vertex). Moreover, the contribution to the latter sum from 
$\bigcup_{k' \le l} E^{k'}_2$ is by
intervals $\widetilde{I}_{k'}$ of
second exploration of some $w_k$, such that $k' \in (k,l]$ and  
\[
F(\cn) - \E(|\widetilde{I}_{k'}| \, |\calF_{k'-1}) 
= H_0 (\widetilde{I}_{k}-X_{k'}) 
\]
measures for an independent interval ${\cal I}' \subseteq \S_{\cn}$ 
of length law $\Gamma_{\cn}(2,1)$ around $X_{k'}$, 
the expected length of all but the part of ${\cal I}'$ 
built within $(\widetilde{I}_{k})^c$ 
(note that for $k' \le N$, a second exploration 
of the vertex requires an a-priori active point 
on it, with $X_{k'} \notin \widetilde{I}_{k}$ then 
$\calF_{k'-1}$-measurable). These precise terms 
appear also in \eqref{eq:bn-expres2} unless
the link $X_{k'}$ has been formed before step $k$, and we
recall Remark \ref{rem:sur-bd} that at most $\Delta_{\rm e}(l)={\rm O}_{\P}(1)$ 
such terms may be missing from \eqref{eq:bn-expres2}. In contrast, the 
$\Delta_{\rm a}(l)$ active points on explored 
vertices after step $l$ do not contribute to the 
left sum of \eqref{eqn:B-N:simplify}, while some 
may participate in \eqref{eq:bn-expres2}. However, as we show in the sequel, 
\begin{equation}\label{eq:bd-Delta-active}
\E [ \Delta_{\rm a} (l) ] \le \frac{l}{N} \E [A_l] \,,\qquad \forall \, l \ge 1 \,,
\end{equation}
hence $\Delta_{\rm a}(l) = {\rm O}_{\P}(1)$ 
uniformly in $l \le s_0 N^{2/3}$ (thanks to Proposition \ref{prop:nu-1-0}).
Thus, as claimed, we have at most
\begin{equation}\label{eq:Deltas}
\Delta_{\rm e}(l)  + \Delta_{\rm a}(l) + 2 \, \nu^l_{\ge 3} = {\rm O}_{\P} (1) \,,
\end{equation}
differing terms between these two sums (recall Proposition \ref{prop:nu-m}
that $\nu^l_{\ge 3} = {\rm O}_{\P}(l^3/N^2) = {\rm O}_{\P}(1)$
uniformly over $l \le s_0 N^{2/3}$). 

Returning to establish \eqref{eq:bd-Delta-active}, note that 
the fraction $q_l$ of explored vertices after step $l$, never
exceeds $l/N$ and the average over vertices of the 
mean number of active points per vertex is $\E[A_l]/N$.
Thus, the bound \eqref{eq:bd-Delta-active} holds 
if these mean numbers are the same across all vertices,  
or more generally, if they are
tilted in favor of the non-explored vertices. 
Further, suffices to 
consider $l \le N$, whereby only parts (a) and (b) of the rules for 
choosing the 
explored points are ever used, see Subsection \ref{subsec:exploration}.
Exploring there the connections of $I_k$,
one uses the same rate to each vertex 
$i \ne w_k$, before erasing some of the connections to explored vertices, 
while erasing none of those to non-explored vertices. 
Beyond this reduction in $\E[\Delta_{\rm a}(l)]$ relative to
$\E[A_l]$, and possibly having $q_l< l/N$ (due to events
$\{ E^k_{\ge 2} \}$), the 
only other deviation from 
uniformity is due to the non-random preference of choosing $w_k$ according 
to its index. The latter can only cause 
the probability that $w_k$ is an already explored vertex (which thereby 
reduces $\Delta_{\rm a}(k)$ by one), to exceed $q_k$. In conclusion, each 
of these effects merely tilts the mean number of active points 
per vertex, towards the non-explored vertices, hence collectively they 
merely reinforce the inequality \eqref{eq:bd-Delta-active}.

The other contribution to 
the ${\rm O}_{\P}(1)$ of \eqref{eq:bn-expres1} comes from the case of 
\begin{equation}\label{eq:bd-E1-E1c}
{\bf 1}_{(E^{k}_1 \cap E^{k'}_1 \cap \{k' \ne k\})^c} \le 
{\bf 1}_{\{k'=k\}} + {\bf 1}_{E^{k}_{\ge 2}} + {\bf 1}_{E^{k'}_{\ge 2}} \,,
\end{equation}
in the right-most sum of \eqref{eqn:B-N:simplify}. Indeed,
uniformly over $l \le s_0 N^{2/3}$,
\begin{align}
\frac{1}{N} \sum_{1 \le k \le k' \le l} ({\bf 1}_{k'=k} 
+ {\bf 1}_{E^{k}_{\ge 2}}
+ {\bf 1}_{E^{k'}_{\ge 2}})
& \le \frac{2l}{N} (1 + \nu^l_{\ge 2}) \nonumber \\
&= \frac{2l}{N} + {\rm O}_{\P}(l^3 N^{-2}) 
=
{\rm O}_{\P}(1)\,. 
\label{eq:bd-Egeq2}
\end{align}
Further, 
the event $E_1^{k} \cap E_1^{k'}$ results for $k'>k$ with $w_{k'} \ne w_{k}$ and
$\widetilde{I}_{k'}$ lying  
on a circle which is completely neutral (other than active points), 
after step $(k'-1)$. So 
following our convention in making $X_{k'}$ an 
$\calF_{k'-1}$-measurable variable, this case contributes  
\[ 
\E ( |\widetilde{I}_{k'} \cap \widetilde{I}_{k}| | \calF_{k'-1} ) 
= H(\widetilde{I}_{k}-X_{k'}) \,,
\]
to the right-most sum of \eqref{eqn:B-N:simplify},
thereby completing the proof.
\end{proof}

\begin{proof}[Proof of Proposition \ref{lemma:B_bar}]
Plugging $F(\c_N) = \lambda_N\left( 1+ a\,N^{-1/3}\right)$, we find 
that the \abbr{rhs} of \eqref{eq:bn-expres1} is 
$$
{\rm O}_{\P}(1) + a l \, N^{-1/3} \,,
$$
which upon setting $l=s N^{2/3}$ and scaling by $N^{-1/3}$, converges 
to $a s$ when $N \to \infty$, uniformly over 
$s \le s_0$. Since $\la_N \to \la_\star = F(\c_\star)$
(see \eqref{eq:window-curve}),
we complete the proof of the proposition
by way of Lemmas \ref{prop:intersection} and \ref{prop:third-sum},
which show that upon 
scaling by $N^{-1/3}$, for $l=s N^{2/3}$ 
the terms in \eqref{eq:bn-expres3} and \eqref{eq:bn-expres2} 
converge in probability as $N \to \infty$, uniformly over $s \le s_0$, 
to the appropriate non-random limits, respectively. Indeed, by 
the preceding we conclude that   
the processes $s \mapsto \bar{B}_N(s)$ converge 
to the deterministic path $s \mapsto \rho^{a,\c_\star}(s)$ 
of \eqref{eqn:def:rho-a-beta}, uniformly over $s \le s_0$.
%
\end{proof}

\begin{lemma}
\label{prop:intersection}
For $(\cn,\la_N)$ satisfying \eqref{eq:window-curve} 
and $l_N=s N^{2/3}$, we have 
\begin{equation}
\label{eqn:intersection2}
\lim_{N\rightarrow\infty} \frac{1}{{l_N \choose 2}} \sum_{k=1}^{l_N}
\sum_{k'=k+1}^{l_N}
H(\widetilde{I}_{k}-X_{k'}) \,
= 
\frac{F^2(\c_\star)}{\c_\star} \,,
\end{equation}
in probability, uniformly over $s \in [\epsilon,s_0]$.
\end{lemma}
\begin{lemma}\label{prop:third-sum}
For $(\cn,\la_N)$ satisfying \eqref{eq:window-curve}, 
$l_N=s N^{2/3}$ 
\begin{align}
 \lim_{N\to\infty} \frac{2N}{l_N^2} 
\sum_{k=1}^{l_N} \sum_{i=1}^{r(k,l_N]}  
H_0(\widetilde{I}_{k}-X_{(k,i)})   
= F(\c_\star) \big(1 - \frac{F(\c_\star)}{\c_\star} - \gamma(\c_\star) \big)\,, 
\label{eqn:third-sum}
\end{align}
in probability, uniformly over $s \in [\epsilon,s_0]$.
\end{lemma}
\begin{remark}\label{rem:heur-inter} In view of Lemma \ref{lem:geo},
the heuristic behind \eqref{eqn:intersection2} is that for 
most pairs $k' > k$ the conditional law of $X_{k'}$ 
given ${\calF}_{k^-}$ is nearly uniform.
To see why we expect \eqref{eqn:third-sum} to hold, let 
\begin{equation}\label{eq:Rl-dfn}
R_l := \sum_{k=1}^l r(k,l] \,,
\end{equation}
be the aggregate number of links by the end of the $l$-th step, 
so $R_k-R_{k-1}$ merely counts the connections made during 
the $k$-th step to all explored vertices, other than $w_k$. As
such, conditional on ${\mathcal F}_{k^-}$, the law of $R_k-R_{k-1}$
is Poisson of rate $|\widetilde{I}_k| q_{k^-}/\lambda_N$, where
$q_{k^-}$ denotes the 
fraction of explored vertices, other than $w_k$, by step $k$.
With $q_{k^-} = k/N (1+o(1))$  
and $F(\cn)/\lambda_N \to 1$, the totality of these Poisson rates
\begin{equation}\label{eq:comp-dfn}
V_l := \frac{1}{\lambda_N} \sum_{k=1}^l |\widetilde{I}_k| \, q_{k^-} \,,
\end{equation} 
should be about $l^2/(2 N)$. We further expect most positions $X_{(k,i)}$ 
to be nearly uniform on $\S^{w_k}_{\cn}$ and approximately independent 
of the first explored interval $\widetilde{I}_{k}$ on that vertex. 
Upon justifying these two approximations, we get \eqref{eqn:third-sum} 
from the \abbr{lln} for the empirical average of 
$\mathsf{UH}_0(\cdot)$
at the nearly i.i.d. $|\widetilde{I}_{k}|$.
\end{remark}

\begin{proof}[Proof of Lemma \ref{prop:intersection}] While proving  
\eqref{eq:UHav}, we saw that at $\cn \to \c_\star$ 
\[
\mathsf{UH}(J) 
= \frac{F(\cn)}{\cn} J 
\to \frac{F(\c_\star)}{\c_\star} J \,.
\]
Thus, by the \abbr{lln}, if $J_{k} := |\widetilde{I}_{k}|$ are i.i.d. 
$\Gamma_{2,1}(\cn)$ variables, then
\begin{equation}\label{eq:conv-G-U-I'} 
\frac{1}{{l_N \choose 2}} \sum_{k=1}^{l_N}
\sum_{k'=k+1}^{l_N}
\mathsf{UH}(J_k) 
= 
\frac{2 F(\cn)}{\cn (l_N-1)} 
\sum_{k=1}^{l_N} (1-\frac{k}{l_N}) J_{k}
\to
\frac{F^2(\c_\star)}{\c_\star} \,,
\end{equation}
in probability, uniformly over $N^{-2/3} l_N \in [\epsilon,s_0]$.
This indeed is the joint law of 
$\{J_{k}\}$ for the \abbr{qrg}, apart from possibly at most
\[
\nu^{l_N}_{\ge 2} = {\rm O}_{\P}(l_N^2/N) = {\rm o}_{\P}(l_N)
\] 
values of $k$. Hence, with $H(\cdot)$ uniformly bounded, the  
uniform convergence in probability of \eqref{eq:conv-G-U-I'}
extends to our setting. Recalling \eqref{eq:UH} and
using the notations ${\bar l}_N := s_0 N^{2/3}$,
\[
\widehat{H}({\cal I}):=H({\cal I})-\E[H({\cal I}-U)|{\cal I}] \,,
\]
with $U \in \S_{\cn}$ uniform, it suffices for proving the lemma, 
to show that 
\begin{equation}\label{eq:2ndmom}
\frac{1}{{\bar l}_N^4} \, \sum_{1 \le k < k' \le {\bar l}_N} \sum_{1 \le j < j' \le 
{\bar l}_N}
\big| \E [\widehat{H}(\widetilde{I}_{k}-X_{k'}) \widehat{H}(\widetilde{I}_{j}-X_{j'}) ] \big| 
\to 0
\end{equation}
Assuming \abbr{wlog} that $j' > r := j \vee k'$ 
(there are only $O({\bar l}_N^3)$ terms with $j'=k'$), and noting that
the uniformly bounded $\widehat{H}(\widetilde{I}_{k}-X_{k'})$ is 
${\calF}_{r^-}$-measurable, it suffices for 
\eqref{eq:2ndmom} to show that  
\begin{equation}\label{eq:h-c}
\frac{1}{{\bar l}_N^2} \, \sum_{1 \le r < r' \le {\bar l}_N} 
h (r,r') \to 0 \,,
\end{equation}
where
\begin{equation}\label{eq:tilH-L1c}
h (r,r') := \max_{j \le r} \Big\{
\E \Big[ \big| \E [\widehat{H}(\widetilde{I}_{j}-X_{r'}) | {\calF}_{r^-} ] \big|
\Big] \Big\} \,.
\end{equation}
Further, 
with $\pi_N$ denoting the uniform law on $\S_{\cn}$ 
we have that for some $C<\infty$ and all $r,r',N$, 
\begin{equation}\label{eq:QN-def}
h(r,r') \le C \,
\E \big[ d_{\rm TV} (\pi_{r,r'},\pi_N) \big] \,,\qquad 
\pi_{r,r'}(\cdot) := \P (X_{r'} \in \cdot | {\calF}_{r^-}) \,.
\end{equation}
In view of \eqref{eq:unif-new} of 
Lemma \ref{lem:lim:unif2},
plugging this bound into \eqref{eq:h-c} concludes the  
proof of Lemma \ref{prop:intersection}. 
\end{proof}
  
\begin{proof}[Proof of Lemma \ref{prop:third-sum}] Recalling that 
if $J_k := |\widetilde{I}_k|$ is $\Gamma_{2,1}(\cn)$ variable
then $\E J_k = F(\cn) = \lambda_N (1+o(1))$, it follows from the \abbr{lln}
by the same argument as in our proof of Lemma \ref{prop:intersection}, 
that the following proxy 
\[
V_l^\star := \frac{1}{\lambda_N} \sum_{k=1}^l J_k \,\frac{k}{N} \,,
\]
for $V_l$ of \eqref{eq:comp-dfn} is such that 
uniformly over $s \in [\epsilon,s_0]$,
\[
\E \Big[ \, \big| \frac{2N}{l_N^2} V^\star_{l_N} - 1 \big| \, \Big] \to 0 \,.
\]
Moreover, for any $k \le l_N$,
\[
\frac{k}{N} - q_{k^-} \in [0,N^{-1} \nu^{l_N}_{\ge 2}] = O_\P(N^{-2/3}) \,.
\]
Hence, $\E [|V_{l_N}^\star - V_{l_N}|] \le O_\P(1)$ and consequently,
uniformly over $s \in [\epsilon,s_0]$,
\begin{equation}\label{eq:qvar-conv}
\E \Big[ \, \big| \frac{2N}{l_N^2} V_{l_N} - 1 \big| \, \Big] \to 0 \,.
\end{equation}
As explained in Remark \ref{rem:heur-inter},  
$M_k:=R_{k-1}-V_{k-1}$ is an ${\mathcal F}_{k^-}$ martingale
and the corresponding predictable part in Doob's decomposition for 
$M_k^2$ is precisely $\langle M \rangle_k = V_{k-1}$.
Thus, by \eqref{eq:qvar-conv}, uniformly over $s$ as above, 
\[
\E \Big[ \big( \frac{2N}{l_N^2} M_{l_N} \big)^2 \Big] 
\to 0 \,,
\]
hence by yet another application of \eqref{eq:qvar-conv}, also  
\begin{equation}\label{eq:Rl-conv}
\E \Big[ \, \big| \frac{2N}{l_N^2} R_{l_N} - 1 \big| \, \Big] \to 0 \,.
\end{equation}
Consider the successive 
steps  
$\varpi_1 \le \varpi_2 \le \ldots \le \varpi_j\le \ldots$ in which links of type 
$\{X_{(k,i)}, \, k,i \ge 1 \}$
to previously explored vertices, say $X^\star_j := X_{(k,i)}$, are formed,
with the induced stopped-filtration $\{{\calF}_{\varpi_j}\}$. We can ignore 
here ties (namely, $\varpi_{j+1} = \varpi_{j}$), since the 
expected number of such is bounded uniformly in $N$ and $l \le s_0 N^{2/3}$
(indeed, $R_{k}-R_{k-1}$ which is at most a 
Poisson of rate $\cn k/(\la_N N)$, yields $(R_k-R_{k-1}-1)_+$ ties).  
Thanks to \eqref{eq:Rl-conv}, it suffices for \eqref{eqn:third-sum} 
to show that for $\ell_N = t N^{1/3}$, 
\begin{align}
 \lim_{N\to\infty} \frac{1}{\ell_N} \sum_{j=1}^{\ell_N}   
H_0(\widetilde{I}_{k(j)}-X^\star_{j})   
= F(\c_\star) \big(1 - \frac{F(\c_\star)}{\c_\star} - \gamma(\c_\star) \big)\,, 
\,\text{ in }\,\P \,,
\label{eqn:av-third-sum}
\end{align}
uniformly over $t \in [\epsilon',t_0]$, where 
$k(j)$ denotes the first 
exploration step of the vertex on which $X_j^\star$ lies. Further, 
for at most $\nu^{\bar{l}_N}_{\ge 3} =O_\P(1)$ of the links $\{X_j^\star, j \le  \ell_N\}$, 
the vertex $w_{k(j)}$
appears more than once in this collection. Consequently, appealing 
to the \abbr{RHS} of \eqref{eq:UHav},
as in the proof of Lemma \ref{prop:intersection} we get analogously to 
\eqref{eq:conv-G-U-I'} that as $N \to \infty$, uniformly over $t$ as above,
\[
\frac{1}{\ell_N} \sum_{j=1}^{\ell_N}
\mathsf{UH}_0(|\widetilde{I}_{k(j)}|) \to
F(\c_\star) \big(1 - \frac{F(\c_\star)}{\c_\star} - \gamma(\c_\star) \big)\,,
\,\text{ in }\,\P \,.
\]
Thus, 
similarly to the proof of Lemma \ref{prop:intersection}, setting
${\bar \ell}_N = t_0 N^{1/3}$ and  
\[
\widehat{H}_0({\cal I}):=H_0({\cal I})-\mathsf{UH}_0(|{\cal I}|) \,,
\]
we get \eqref{eqn:av-third-sum} as soon as we show that  
\begin{equation}\label{eq:2ndmom-part2}
\frac{1}{{\bar \ell}_N^2} \, \sum_{1 \le j < j' \le {\bar \ell}_N} 
\big| \E [\widehat{H}_0(\widetilde{I}_{k(j)}-X^\star_{j}) \widehat{H}_0(\widetilde{I}_{k(j')}-X^\star_{j'}) ] \big| 
\to 0 \,.
\end{equation}
With  
$\widehat{H}_0(\widetilde{I}_{k(j)}-X^\star_{j})$ 
uniformly bounded and  
${\calF}_{\varpi_j}$-measurable, this in turn follows from 
\begin{equation}\label{eq:h0-c}
\frac{1}{{\bar \ell}_N^2} \, \sum_{j=1}^{{\bar \ell}_N} \sum_{j'=j+1}^{{\bar \ell}_N} 
\E \Big[ \big| \E [\widehat{H}_0(\widetilde{I}_{k(j')}-X^\star_{j'}) | 
{\calF}_{\varpi_j} ] \big|
\Big] 
\to 0 \,.
\end{equation}
Further, setting $r=r(j,j'):=\varpi_j \vee k(j')$,
\begin{equation}\label{eq:Q0N-def}
\pi^\star_{j,j'}(\cdot) := \P (X^\star_{j'} \in \cdot | \calF_{r(j,j')}) \,,
\end{equation}
and recalling the definition \eqref{eq:UH} of $\mathsf{UH}_0(\cdot)$
with $U \sim
\pi_N$ on $\S_{\cn}$,  
we have that for some $C<\infty$ and all $j,j',N$, 
\[
\E \Big[ \big| \E [\widehat{H}_0(\widetilde{I}_{k(j')}-X^\star_{j'}) | 
{\calF}_{\varpi_j} ] \big|
\Big] 
\le C \,
\E \big[ d_{\rm TV} (\pi^\star_{j,j'},\pi_N) \big] \,.
\]
In view of 
Lemma \ref{lem:lim:unif3},
plugging this 
into \eqref{eq:h0-c} concludes the proof.
\end{proof}

We next complete the proof of Lemma \ref{prop:intersection} by showing 
that for most $r<r' \le {\bar l}_N$, the law of $X_{r'}$ given $\calF_{r^-}$
is nearly uniform (in total-variation distance).
\begin{lemma}
\label{lem:lim:unif2}
Setting ${\bar l}_N = s_0 N^{2/3}$ we have for 
$\pi_{r,r'}$ of \eqref{eq:QN-def}, that 
\begin{align}
\label{eq:unif-new}
\lim_{N \to \infty} 
\frac{1}{{\bar l}_N^2}
\sum_{r = 1}^{{\bar l}_N} \sum_{r'=r+1}^{{\bar l}_N} 
\E \Big[ d_{\rm TV} (\pi_{r,r'},\pi_N) \Big] 
= 0 \,.
\end{align}
\end{lemma}
\begin{proof} Fixing $r \ge 1$
%
we assign to each active point at the end of the $r$-th step, and to
any link formed thereafter (including a link onto the previously 
explored part of the \abbr{qrg}), a proxy counter for its 
uniformity and ${\calF}_{r}$-independence, as follows. First, 
the counter of each active point at the end 
of the $r$-th step is set to $0$ if it is on an unexplored vertex 
and to $-1$ otherwise. Then, sequentially in $k>r$, 
if $A_{k-1}>0$ and the vertex $w_k$ was not previously explored, 
we set as $n_{r,k}$ the counter of the 
active point $(w_k,X_k)$, if $A_{k-1}=0$ we set 
$n_{r,k}=\infty$, and otherwise let $n_{r,k}=-1$. Thereafter,
each link formed during the $k$-th step gets a counter
value $n_{r,k}+1$. We claim that for 
$\rho:=\sup_N \{ \P(\Gamma_{\cn} (2,1) < \cn) \} < 1$ 
and any $r' \in (r,N)$, 
\begin{equation}\label{eq:counter-new}
d_{\rm TV} (\P (X_{r'} \in \cdot | {\calF}_{r^-}, n_{r,r'} \ge c),\pi_N) 
\P(n_{r,r'} \ge c | {\calF}_{r^-} ) \le \rho^{c} \,.
\end{equation}
Indeed, at each step we choose the value of $X_k$ independently of the 
positions of active points within their respective circles. Having
$A_{k-1}=0$ yields $X_k \sim \pi_N$ independently of ${\calF}_{k-1}$
and this property is inherited by any point 
to which the path from $X_k$ 
involves only first explorations of the relevant vertices. Further, 
the ${\calF}_{k-1}$-measurable event $n_{r,k} \ge 0$ (namely $E^k_1$),
results with
$\P(\widetilde{I}_k=\S_{\cn}|{\calF}_{k-1}) \ge 1-\rho$ and
$\widetilde{I}_k=\S_{\cn}$ yields by ${\cal L}_{w_k,i}$  
uniformly distributed links to all
vertices (prior to the erasures on previously explored space). 
In case $n_{r,r'} \ge c$ is finite, one has during $(r+1,r')$ 
at least $c$ consecutive forefathers 
of $(w_{r'},X_{r'})$ in our exploration tree, all of whom were
explored on circles which are neutral (apart from active points). 
The chance that 
none of these forefathers forced a uniform conditional 
law of $X_{r'}$, is by the preceding at most $\rho^c$, 
thereby establishing \eqref{eq:counter-new}. 

Now, thanks to \eqref{eq:counter-new} and the convexity of
the $[0,1]$-valued $d_{\rm TV} (\cdot,\pi_N)$, it suffices 
for \eqref{eq:unif-new} to show that for any $c<\infty$,
\begin{equation}\label{eq:high-new-counter}
\lim_{N \to \infty} \E [ \Gamma_c(\bar l_N) ] = 0 \,, \quad \textrm{where} \quad
\Gamma_c(l) := l^{-2} \sum_{1 \le r < r' \le l} \, {\bf 1}_{\{n_{r,r'}  < c\}}  \,.
\end{equation}
To this end, with $\calZ_c(k')$ denoting the size of the exploration 
sub-tree of depth at most $c$, rooted at the active point $X_{k'}$
and $\Theta(r)$ enumerating those $k' > r$ for which the link to $(w_{k'},X_{k'})$ 
has been formed before the end of the $r$-th step, we claim that
\begin{equation}\label{bd:Gamma-c}
\Gamma_c(l) \le 
l^{-2} \sum_{1 \le r < k' \le l} \calZ_c(k') {\bf 1}_{\{k' \in \Theta(r)\}} +
l^{-1} \sum_{k'=1}^l \calZ_c(k') {\bf 1}_{E^{k'}_{\ge 2}} 
\,.
\end{equation}
Indeed, to have $n_{r,r'} < c$, one of the $c$ consecutive 
exploration forefathers of $(w_{r'},X_{r'})$, say $(w_{k'},X_{k'})$,
must have been an event $E^{k'}_{\ge 2}$ (namely, not a first exploration), 
or alternatively, be an active point formed before the end of the $r$-th step.
In the latter case, considering the \emph{last such step} 
(i.e. along the path to $(w_{r'},X_{r'})$ and among active points 
formed by the end of the $r$-th step), guarantees that $r<k'$ 
in the first sum 
of \eqref{bd:Gamma-c}. Thereafter, 
$\calZ_c(k')$ bounds the number of possible pairs $(k',r')$ having
path distance at most $c$
and should no such $(k',r')$ exist, the second sum on the \abbr{rhs} of \eqref{bd:Gamma-c} bounds the number of $r' \le l$ with some
previously explored vertex among the last $c$ steps on the exploration 
path to $(w_{r'},X_{r'})$.

Now, for any $k'$, conditional on $\calF_{k'^{-}}$ the 
variable $\calZ_c(k')$ is stochastically dominated by the size of
a Galton-Watson tree of depth $c$ and a Poisson($\mu$) off-spring law, for 
$\mu := \max_N \{\c_N/\la_N\}$ finite. It thus follows that for some 
$\kappa_c=\kappa_c(\mu)$ finite and all $N$,
\[
\max_{k' \ge 1} \E [\calZ_c(k')|\calF_{k'^{-}} ] \le \kappa_c \,.
\]
Equipped with the latter bound, upon 
considering the expected values in \eqref{bd:Gamma-c}, since
both $E^{k'}_{\ge 2}$ and $\{k' \in \Theta(r)\}$ are in $\calF_{k'^{-}}$, 
while $|\Theta(r)| \le A_r$, it follows by the tower property of the 
conditional expectation, that 
\begin{equation}\label{bd:exp-Gamma-c}
\E [\,\Gamma_c(l)\,] \le \frac{\kappa_c}{l} \Big( 
\E \big[ \max_{r < l} \{A_r\} \big] +
\E \big[ \sum_{k'=1}^l {\bf 1}_{E^{k'}_{\ge 2}} \big]
\Big) \,. 
\end{equation}
Further, recall that $\sum_{k' \le l} {\bf 1}_{E^{k'}_{\ge 2}} \le \nu^l_{\ge 2}$. Thus,
with both $l^{-1} \E [ \, \max_{r<l} A_r \, ]$ and
$l^{-1} \E [ \, \nu_{\ge 2}^{l} \, ]$ decaying to zero at
$l={\bar l}_N$ and $N \to \infty$ (due to 
Proposition \ref{prop:nu-1-0} and
\eqref{eq:num-exp}, respectively),
the bound \eqref{bd:exp-Gamma-c} yields that \eqref{eq:high-new-counter} holds.
\end{proof}

Similarly to Lemma \ref{lem:lim:unif2}, we complete the proof of Lemma \ref{prop:third-sum} by showing 
that for most $j<j' \le {\bar \ell}_N$, the conditional law of $X^\star_{j'}$
is nearly uniform.
\begin{lemma}
\label{lem:lim:unif3}
Setting ${\bar \ell}_N = t_0 N^{1/3}$ we have for $\pi^\star_{j,j'}$ 
of \eqref{eq:Q0N-def}, that 
\begin{align}
\lim_{N \to \infty} 
\frac{1}{{\bar \ell}_N^2}
\sum_{j = 1}^{{\bar \ell}_N} \sum_{j'=j+1}^{{\bar \ell}_N} 
\E \Big[ d_{\rm TV} (\pi^\star_{j,j'},\pi_N) \Big] = 0 \,.
\label{eq:unif-link}
\end{align}
\end{lemma}
\begin{proof} We record the 
first exploration step $\mathsf{ex}(w) \ge 1$ of each vertex $w \in [1,N]$, 
so upon forming a link of type
$X^\star_{j'}$ onto a vertex $v_{j'}$ one has
that $k(j')=\mathsf{ex}(v_{j'})$. Utilizing the proxy counters of
Lemma \ref{lem:lim:unif2}, let $n^\star_{j,j'}$
denote the value of the counter for the link $X^\star_{j'}$ starting at step 
\[
r = r(j,j') := \varpi_j \vee \mathsf{ex}(v_{j'}) \,. 
\]
By the same reasoning as in the derivation of \eqref{eq:counter-new}, we have that
\begin{equation}\label{eq:counter-link}
d_{\rm TV} (\P (X^\star_{j'} \in \cdot | {\calF}_{r}, n^\star_{j,j'} \ge c),\pi_N) 
\P(n^\star_{j,j'} \ge c | {\calF}_{r} ) \le \rho^{c} \,.
\end{equation}
We thus get \eqref{eq:unif-link} by establishing the 
analog of \eqref{eq:high-new-counter}. That is, upon showing that 
for $l := 2 \sqrt{\ell N}$ and any $c<\infty$,  
\begin{equation}\label{eq:high-link-counter}
\lim_{N \to \infty} \E[\Gamma^\star_c({\bar \ell}_N)] = 0  
 \,, \quad \textrm{for} \quad
\Gamma^\star_c(\ell) := \ell^{-2} {\bf 1}_{\{\varpi_\ell \le l\}} 
\sum_{1 \le j < j' \le \ell} \, {\bf 1}_{\{n^\star_{j,j'}  < c\}} 
\end{equation}
(by \eqref{eq:Rl-conv} it suffices to consider 
$2 N \ell /\varpi_\ell^2 \to 1$, hence the restriction here to 
$\varpi_\ell \le 2 \sqrt{\ell N}$).
Next, for $k < k'$, let $\calZ_c(k';k)=\calZ_c(k')$ if the exploration
sub-tree of depth at most $c$ rooted at the active point $X_{k'}$, 
has a link to $\S_{\cn}^{w_k}$, otherwise setting
$\calZ_c(k';k) \equiv 0$. 
Likewise, $\calZ^\star_c(k'):=\max_{k < k'} \{\calZ_c(k';k)\}$ is the 
non-zero $\calZ_c(k')$ iff the relevant sub-tree produces a link to some 
previously explored vertex.
Setting hereafter $k=k(j')$, recall 
that $n_{r,r'} < c$ requires that one of the $c$ consecutive 
exploration forefathers of the link $X^\star_{j'}$, say $(w_{k'},X_{k'})$,
must have been an event $E^{k'}_{\ge 2}$ (namely, not a first exploration), 
or alternatively, be an active point formed before the end of the 
$r(j,j')$-th step. We thus claim, 
similarly to \eqref{bd:Gamma-c}, that 
\begin{align}\label{bd:Gamma-star-c}
\Gamma^\star_c(\ell) \le 
\ell^{-2} \sum_{j < \ell, \varpi_j <  k' \le l} \calZ^\star_c(k') 
{\bf 1}_{\{k' \in \Theta(\varpi_j)\}} 
&+
\ell^{-1} \sum_{1 \le k < k' \le l} \calZ_c(k';k) 
{\bf 1}_{\{k' \in \Theta(k)\}}  
\nonumber \\
&+ \ell^{-1} \sum_{k'=1}^l \calZ^\star_c(k') {\bf 1}_{E^{k'}_{\ge 2}} \,. 
\end{align}
The first two expressions on the \abbr{rhs} 
distinguish having $k < \varpi_j$ from the case of $k \ge \varpi_j$,
where we sum over $j \le \ell$ 
and cover all choices of $j'$ by the additional sum over $k=k(j') < l$. 
As done on the \abbr{rhs} of \eqref{bd:Gamma-c}, in both expressions we guarantee
that $r<k'$ by having $(w_{k'},X_{k'})$ stand for the last active point on 
the path to $X^\star_{j'}$ among those formed by the end of the $r$-th step. 
The link $X^\star_{j'}$ must lie on the vertex $w_{k(j)}$, yielding the bound 
$\calZ_c(k';k)$ in case $k=r<k'$ with an exploration path distance at most $c$ from
$X_{k'}$ to $X^\star_{j'}$. However, in case $k < \varpi_j =r$ we do not keep 
track of $k$, hence must replace $\calZ_c(k';k)$ by the larger $\calZ^\star_c(k')$ 
which only indicates the existence of a point 
of exploration path distance at most $c$ from $X_{k'}$ which is on a vertex 
that was first explored prior to step $k'$. Finally, should no active point 
$(w_{k'},X_{k'})$ with $k' > k$ of exploration path distance at most $c$ 
from $X^\star_{j'}$ be formed by the end of the $r$-th step, the  
last sum on the \abbr{rhs} of \eqref{bd:Gamma-star-c}
bounds the number of $j' \le \ell$ (with $\varpi_{j'} \le l$),
having a non-neutral circle (explored at some step $k' > k$),   
among the last $c$ points on the 
path to $X^\star_{j'}$. 

Next, recall that $\calZ_c(k')$ is, conditionally on $\calF_{k'^{-}}$,
stochastically dominated by the size of a depth $c$
Galton-Watson tree of a Poisson($\mu$) off-spring law.
The latter size variable has finite moments of all order, whereas $\calZ_c(k';k)$ 
further demands having at least one tree vertex corresponding 
to the prescribed $w_k$. With the production of the specified vertex $w_k$
stochastically dominated by a Poisson of  
rate $\mu/N$, we have for some $\kappa_c^\star$ finite and all $N$,
\[
\max_{k < k'} \{
\E[ \,\calZ_c (k';k)\,|\,\calF_{k'^{-}} \,] \} \le \frac{\kappa^\star_c}{N}\,, \qquad 
\max_{k' \le l} \{ \E[ \,\calZ_c^\star(k')\,|\,\calF_{k'^{-}} \,] \} \le \frac{\kappa^\star_c l}{N}
\,.
\]
Now, analogously to the derivation of \eqref{bd:exp-Gamma-c}, upon 
taking the expectation on both sides of \eqref{bd:Gamma-star-c}, we
get by the tower property and the preceding estimates that
\begin{align*}
\E [ \Gamma^\star_c(\ell) ] \le \frac{\kappa^\star_c l}{N \ell} 
\Big( \E\big[ \ell^{-1} \sum_{j<\ell} \, A_{\varpi_j}
{\bf 1}_{\{\varpi_j \le l\}}  \, \big] + 
\E \big[ l^{-1} \sum_{k'<l} A_{k'} \big]  +
\E \big[ \sum_{k'=1}^l {\bf 1}_{E^{k'}_{\ge 2}} \big] \Big) \,. 
\end{align*}
Utilizing the fact that $l/(N \ell) = 4/l$, we arrive at the same bound 
as in the \abbr{rhs} of \eqref{bd:exp-Gamma-c}. Setting 
$\ell=\bar \ell_N$ corresponds to having $l=\bar l_N$, thus yielding 
\eqref{eq:high-link-counter}
by the reasoning provided at the end of the proof of Lemma \ref{lem:lim:unif2}.
\end{proof}


\subsection{Joint convergence of component sizes}
\label{subsec:joint:convergence}

Recall the statement of Theorem \ref{thm:joint:cgs}. In this section, we shall conclude the proof of this 
theorem using results from previous sections.

As pointed out in \cite{aldous}, Theorem \ref{thm:joint:cgs} primarily has two parts:
\begin{itemize}
\item[1.] First is to prove that the excursions of the limit process are matched by the excursions of
the breadth first random walk.
\item[2.] Second is to arrange these excursions in the decreasing order. This can be achieved if one can ascertain that there exists a random
point after which one is sure (with high probability) not to see large excursions. 
\end{itemize}

In order to settle the first issue, we shall invoke 
\cite[Lemmas 7 and 8]{aldous}, which can be applied verbatim to our case,
together with Proposition \ref{thm:scaling:limit} proved in a previous subsection.

Thus, we only need to be concerned about the second issue, for which we shall need to prove an appropriate version of 
\cite[Lemma 9]{aldous} suited to our case.

Like in \cite{aldous}, let us define 
$$T(y) = \min\{s: W^{a,\c_\star}(s) = -y\},$$
$$T_N(y) = \min\{i: Y_N(i) = -\lfloor yn^{1/3}\rfloor\}.$$
Notice that as a consequence of Proposition \ref{thm:scaling:limit}
\begin{eqnarray*}
N^{-2/3} T_N(y) 
& \rightarrow_{d} & T(y).
\end{eqnarray*}

Therefore, the following lemma completes the proof.

\begin{lemma}
\label{lem:9:aldous}
Let us denote by $p(N,y,\delta)$ the probability that the \abbr{qrg} with the parameters $(\c_N,\lambda_N)$ that satisfy 
\eqref{eq:window-curve}, contains a component 
of size at least $\delta N^{2/3}$ which does not contain any vertex $i$ with $1\le i\le y N^{1/3}$. Then,
$$\lim_{y\to\infty} \limsup_{N \to \infty} p(N,y,\delta) = 0 \,\,\,\,\,\,\, \text{ for all } \delta >0.$$
\end{lemma}
\begin{proof} Fix $\delta > 0$. Let $v_{\calC_i}$ be the minimal vertex of the 
component $\calC_i$, then for an interval $\mathfrak{I}\subset\real_+$, define
\begin{equation}
\label{eqn:def:qN}
q(N,\mathfrak{I}) = \E \Big(\sum_{i \ge 1} {\bf 1}_{\left(|\calC_i| \ge \delta N^{2/3}; v_{\calC_i}\in N^{1/3}\mathfrak{I}\right)}\Big).
\end{equation}
Conditioned on arranging the components in a decreasing order of their sizes,
the labels of the vertices of any given fixed 
component $\calC_i$
are going to be uniformly randomly ordered. 
Given such components ordering, define
$$\chi_N(\calC_i) = N^{-1/3} v_{\calC_i},$$
and $${\cal U}_{\calC_i} = N^{-2/3} \left( \text{number of vertices in the component } \calC_i\right)\,.$$
Then note that for any $x \ge 0$
$$\mathbb{P}(v_{\calC_i}>N^{1/3}x \, \big| \, \calU_{\calC_i})=\Big(
1-\frac{{\calU}_{\calC_i}N^{2/3}}{N}\Big)^{N^{1/3}x},$$ 
implying
\begin{equation}
\label{eqn:exp}
\mathbb{P}\left(\chi_N(\calC_i)>y\big|\,{\calU}_{\calC_i}\right) \le
\frac{{\rm e}^{-{\calU}_{\calC_i}y}}{1-{\rm e}^{-{\calU}_{\calC_i}}}\mathbb{P}\left(\chi_N(\calC_i)\leq 1\big|\, {\calU}_{\calC_i}\right) \,.
\end{equation}
Further
\begin{equation}
\P{\left(v_{\calC_i}\in [yN^{1/3},\infty)\right)}=\E\left(\P{\left(v_{\calC_i}\in [yN^{1/3},\infty) \big| \,{\calU}_{\calC_i}\right)}\right).
\end{equation}
At this point, conditional on component sizes being $|\calC_i| = b N^{2/3}$, we note that one can adopt the proof of Proposition 
\ref{prop:nu-m} to the original exploration process restricted to the construction of $\calC_i$ in order to derive similar results 
for $\nu^{\calC_i}_{\ge 3}$, the number of explored intervals belonging to $\calC_i$ sampled by the end of the construction 
of $\calC_i$, which belong to vertices (circles), having at least three such intervals each.
Then, observing that ${\cal U}_{\calC_i}N^{2/3}\geq \frac12(bN^{2/3}-\nu^{\calC_i}_{\ge 3})$, we have
for $\varepsilon>0$,
$$\P\Big({\calU}_{\calC_i}\geq \frac12 b-\frac12\frac{N^{1/3+\varepsilon}}{N^{2/3}}\Big)
\ge 
\P\left(\nu^{\calC_i}_{\ge 3}\leq N^{1/3+\varepsilon}\right)=1-{\rm o}(N^{-1/3}),$$
implying that ${\cal U}_{\calC_i}\in \left(\frac{b}3, b\right)$ with probability $\left( 1-{\rm o}(N^{-1/3}) \right)$.
Consequently,
\begin{eqnarray*}
\E\left(\P\left(\left. v_{\calC_i}\in [yN^{1/3},\infty) \right\vert {\cal U}_{\calC_i}\right)\right) 
& \le & \frac{{\rm e}^{-by/3}}{1-{\rm e}^{-b/3}}\P\left( v_{\calC_i}\in [0,N^{1/3}]\right) +{\rm o}(N^{-1/3}).
\end{eqnarray*}
Recalling the definition of $q(N,\mathfrak{I})$ from \eqref{eqn:def:qN}, and conditioning on the number $M_b$ of components of size $bN^{2/3}$, while observing that
given the sizes of components the minimal vertices of various different components are identically distributed, we get
\begin{eqnarray*}
q(N,[y,\infty))
& = &\E\Big[\sum\limits_{b=\delta}^{\infty}M_b \P\big(\, v_{\calC_i}\in [yN^{1/3},\infty) \,\big\vert \, |\calC_i| = b N^{2/3} \big)\, \Big]\,.
\end{eqnarray*}
With $\sum_{b \ge \delta} M_b \le \delta^{-1} N^{1/3}$, by the preceding, this and \eqref{eqn:exp} imply that
\begin{eqnarray*}
q(N,[y,\infty))
& \le & \frac{{\rm e}^{-\delta y/3}}{1-{\rm e}^{-\delta/3}}q(N, [0,1])+{\rm o}(1).
\end{eqnarray*}
Since $p(N,y,\delta) \le q(N,[y,\infty))$, 
to prove the theorem, it suffices to show that
\begin{equation}
\label{eqn:sup:1}
\sup_N q(N,[0,1]) < \infty,
\end{equation}
Writing $t_i(v)$ as points on the $v$-th vertex around which intervals are constructed and explored, and denoting $\calN(v)$ as the number of such points 
we observe that
$$q(N,[0,1]) \le \sum_{v=1}^{N^{1/3}} \E \Big( \sum_{i=1}^{\calN(v)} {\bf 1}_{\{|\calC(t_i(v))| > \delta N^{2/3}\}}\Big),$$
where $\calC(t_i(v))$ is the maximal connected component containing $t_i(v)$.

Clearly the collection $\{t_1(v),\ldots,t_{\calN(v)}(v)\}$ is independent and identically distributed for different $v\in G_N$.
We replace the exploration by the overcounting process of Section \ref{subsec:overcounting} which is coupled with the exploration
process until the exploration process hits zero. Then, we restart an independent (and identical) overcounting process together with 
restarting the exploration process. We repeat this process until the end of exploration of the complete graph. 
Subsequently, setting $\{t^*_1(v),\ldots,t^*_{\calN^*(v)}\},\,\,\calN^*(v)$ and $\calC^*(t^*_i(v))$ as the corresponding elements of the overcounting process,
we observe that since $|\calC(t^*_i(v))|$ are i.i.d. we have
\begin{eqnarray*}
q(N,[0,1]) &\le & \sum_{v=1}^{N^{1/3}} \E\Big(
\sum_{i=1}^{\calN^*(v)} {\bf 1}_{\{|\calC(t^*_i(v))| > \delta N^{2/3}\}}
\Big)\\
&=& N^{1/3} \E(\calN^*(v)) \P\left( |\calC(t^*_i(v))| > \delta N^{2/3}\right),
\end{eqnarray*}
where we have used Wald's equality.

Therefore, it suffices to prove that $N^{1/3} \P(|\calC^*_0|\geq \delta N^{2/3})$ is bounded by a constant 
where $\calC^*_0$ is a {\it typical} component of the overcounting process.

We now define the coupled overcounting process via
i.i.d. random variables $\xi^w_k$, where each $\xi^w_k$ represents the number of links generated at $k$-th time step by the
overcounting process with the parameters $\lambda_N$ and $\c_N$ lying in the critical window \eqref{eq:window-curve},
Unlike Section \ref{subsec:overcounting}, here
 $\xi^w_i \sim \text{Poisson}\left( \frac{\c_N}{\lambda_N}\right)$.
Then define $S^w_k = S^w_{k-1} + (\xi^w_k - 1)$, with $S^w_0=1$. Setting 
$\tau^w = \min\{k\ge 1: S^w_k = 0\}$, it suffices
to show that $N^{1/3}\P\left( \tau^w > \delta N^{2/3}\right)$ is bounded by a universal constant.
Using same arguments as used in Section \ref{subsec:overcounting}, we 
conclude that
\begin{equation}\label{eqn:overcount-w1}
\P\left( \tau = n+1\right) \,\, \le \,\, \frac{\E(S^w_1)}n \,\,\, \sup_{\ell} \{\P(S^w_{n+1} - S^w_1 = -\ell)\}
\end{equation}

Using \cite[Proposition 2.4.4]{LawlerLimicBook} observe that
$$\P(S^w_{n+1} - S^w_1 = -\ell) \le \frac{c}{n^{1/2}}.$$
Therefore,
$$\P\left( \tau = n+1\right) \,\, \le \,\,c\, n^{-3/2}.$$
Subsequently, following the same arguments as in Section \ref{subsec:overcounting}, we conclude that
$$N^{1/3} \P(|\calC^*_0|\geq \delta N^{2/3})\le c,$$
for some finite $c=c(\delta)$,
thereby proving the statement of the lemma.
\end{proof}

\section{Appendix: Proof of Proposition \ref{prop:nu-1-0}}

With $Y_N=Y_N^{a,\c_N}$, $Z_N=Z_N^{a,\c_N}$ and writing
$\iota(l)$ for the number of maximal connected components in the corresponding graph 
completely explored before step $l$, we use the relations
\begin{align}
A_l  &= Y_N(l) + \iota(l) = Z_N(l) - \textrm{sur}(l) + \iota(l) \,,\\
\iota(l) &= 1- \min_{0\le k\le (l-1)}\{Z_N(k) - \textrm{sur}(k)\} \,, 
\end{align} 
and the fact that $k \mapsto \textrm{sur}(k)$ is non-decreasing, to find that 
\begin{eqnarray*}
A_l 
& = & 1 + Z_N(l) - \textrm{sur}(l) + \max_{k\le (l-1)} \{ \textrm{sur}(k) - Z_N(k) \} \\
& \le &  1 +  \max_{k\le l} \{ Z_N(l) - Z_N(k) \},
\end{eqnarray*}
which can further be simplified to write
\begin{eqnarray*}
A_l &\le & 1 + 2 \max_{k\le l}|Z_N(k)|
\end{eqnarray*}

Recall the martingale decomposition, $Z_N(k) = M_N(k)+ B_N(k)$, where $M_N$ is a martingale and $B_N$ is the
predictable process. Then, for any fixed positive $K$, set
$$\Upsilon_N = \min\{k: |Z_N(k)| > KN^{1/3}\}\wedge (s\,N^{2/3}).$$
By Markov's inequality, it thus suffices for 
Proposition \ref{prop:nu-1-0} to show that
\begin{equation}\label{appendix:ntp-2}
\E \big[ |Z_N(\Upsilon_N)|^2 \big] = {\rm O}(N^{2/3}) \,.
\end{equation}
To this end, using the notation introduced in \eqref{eq:dec12}, clearly
\begin{equation}\label{appendix:mart-decompo}
\E\big[ |Z_N(\Upsilon_N)|^2 \big] \le 2 \E \big[ |M_N(\Upsilon_N)|^2 \big] + 
2 \E \big[ |B_N(\Upsilon_N)|^2 \big].
\end{equation}
Further, by Doob's optional sampling theorem,
$$\E \big(M_N(\Upsilon_N)^2\big) = \E(D_N(\Upsilon_N)) \le \E(D_N(sN^{2/3})),$$
since $M_N^2$ is a sub-martingale and $\Upsilon_N \le sN^{2/3}$.
Now, recall that
\begin{eqnarray*}
D_N(sN^{2/3}) = \sum_{k=1}^{sN^{2/3}} \text{var}(\left.\zeta_k \right| \calF_{k-1})
\end{eqnarray*}
with by way of \eqref{eq:zeta-to-phi} for uniformly bounded
$\varphi_k \le \cn/\la_N$, 
has expected value
bounded by $cN^{2/3}$, for some finite $c=c(a,s)$. 
Turning to show the same for 
$\E |B_N(\Upsilon_N)|^2$, recall that the right sum 
of \eqref{eqn:B-N:simplify} has
$l(l+1)/2$ terms, each bounded by $\cn/N$, whereas to 
the left sum only the at most 
$\nu^l_{\ge 2}$ events $E^{k'}_{\ge 2}$ contribute
(no more than $2\cn$ each). Thus, in view of \eqref{eq:window-curve} and
\eqref{eqn:B-N:simplify},
\begin{equation}\label{eq:last-exp}
B_N(l) \le 1 + \frac{a l}{N^{1/3}} + 
\frac{\cn l^2}{\la_N N} +
 \frac{2 \cn}{\la_N} \nu^l_{\ge 2}   
 \,.
\end{equation}    
For $l= \Upsilon_N \le s N^{2/3}$ the non-random part of the 
\abbr{rhs} of \eqref{eq:last-exp} is at most $c N^{1/3}$. Next, 
upon examining the argument leading to 
\eqref{eq:num-exp}, we deduce that 
$\nu^l_{\ge 2}$ is stochastically dominated by the sum of 
at most $\nu^l_{\ge 1}$ i.i.d. Poisson variables of rate
$\kappa' l N^{-1}$ each. Hence, for some $C$ finite and all $l,N$,
\[
\E \Big[ (\nu^l_{\ge 2})^2 \Big] \le C  + C \frac{l^2}{N^2} 
\E \Big[ (\nu^l_{\ge 1})^2 \Big] \,.
\]
Similar refinement in the argument leading to \eqref{eq:nu1-exp}, yields that
\[
\E \Big[ (\nu^l_{\ge 1})^2 \Big] \le C l^2 \,, 
\]
hence $\E \Big[ (\nu^l_{\ge 2})^2 \Big] \le c N^{2/3}$ for $l \le s N^{2/3}$.
Such bound holds for $\E |B_N(\Upsilon_N)|^2$ and  
the decomposition \eqref{appendix:mart-decompo} 
yields \eqref{appendix:ntp-2} (thereby completing the proof). \qed


\end{document}